\documentclass[10pt,fleqn]{article}
\usepackage{amsmath,amssymb,amsfonts,amsthm,color,dsfont,graphicx}
\usepackage[a4paper,margin=3cm]{geometry}
\usepackage[english]{babel}
\usepackage[utf8]{inputenc}
\usepackage[pdftex]{hyperref}

\graphicspath{{./figures/}}

\allowdisplaybreaks

\newcommand{\R}{\mathds{R}}

\newcommand{\1}{\mathds{1}}

\newcommand{\tdiff}[1]{\frac{\mathsf{d}}{\mathsf{d} #1}}
\newcommand{\dx}{\;\mathsf{d\mathbf{x}}}

\newcommand{\dSx}{\;\mathrm{da(\mathbf{x})}}

\newcommand{\D}{\ensuremath{\mathrm{D}}}

\newcommand{\eps}[1]{\varepsilon[#1]}

\newcommand{\workdom}{D}
\newcommand{\body}{\Omega}
\newcommand{\etensor}{C}
\newcommand{\efftensor}{C^\ast}
\newcommand{\strain}{u}
\newcommand{\stress}{\sigma}

\newcommand{\intga}{\int_{\Gamma_N}\,}

\newcommand{\Vh}{\mathcal{V}^{(2)}_h}
\newcommand{\Ih}[1]{\mathcal{I}^{(#1)}_h}

\newcommand{\Lag}{\mathcal{L}}

\newcommand{\ie}{i.\,e.}
\newcommand{\eg}{e.\,g.}
\newcommand{\wrt}{w.\,r.\,t.}

\newcommand{\fr}[2]{{\textstyle\frac{#1}{#2}}}

\renewcommand{\div}{\mathrm{div}}

\newtheorem{theorem}{Theorem}[section]
\newtheorem{corollary}{Corollary}
\newtheorem{lemma}[theorem]{Lemma}

\theoremstyle{definition}

\newtheorem{remark}{Remark}

\DeclareMathOperator{\tr}{tr}

\begin{document}
\title{A posteriori error estimates for sequential laminates in shape optimization}

\author{
       Benedict Geihe\footnotemark[1]
  \and Martin Rumpf\footnotemark[1]
}

\renewcommand{\thefootnote}{\fnsymbol{footnote}}
\footnotetext[1]{Institute for Numerical Simulation, Rheinische Friedrich-Wilhelms-Universit\"at Bonn, Endenicher Allee 60, 53115 Bonn, Germany}
\renewcommand{\thefootnote}{\arabic{footnote}}

\maketitle

\begin{abstract}
A posteriori error estimates are derived in the context of two-dimensional structural elastic shape optimization under the compliance objective.
It is known that the optimal shape features are microstructures that can be constructed using sequential lamination.
The descriptive parameters explicitly depend on the stress.
To derive error estimates the dual weighted residual approach for control problems in PDE constrained optimization is employed, 
involving the elastic solution and the microstructure parameters.
Rigorous estimation of interpolation errors ensures robustness of the estimates while local approximations are used to obtain 
fully practical error indicators. Numerical results show sharply resolved interfaces between regions of full and intermediate material density.\\
\end{abstract}

\section{Introduction}\label{sec:intro}
In this paper we will address error estimation for shape optimization problems and study an adaptive finite element method based on
an optimal microstructured elastic material composed of sequential laminates. A classical shape optimization problem amounts to 
distribute a fixed amount of rigid material throughout a domain in order to optimize a particular cost functional
when confronted with surface and volume loads. However it is known that the problem is ill-posed if only scale invariant cost functionals
are taken into account. As the formation of microstructures comes along with this deficiency most approaches 
enforce regularity by an additional regularizing cost functional such as the shape perimeter. Alternatively the problem
can be relaxed to allow for more general shapes characterized by local effective material properties and densities. They stem from
the shape structure on the microscale influencing the effective macroscopic description via homogenization. Optimal microstructures exist but are not unique.
A particular flexible microstructure is constructed via sequential lamination. This allows for a convenient algorithmic treatment as a small set of parameters
characterizing the local microstructure can be computed explicitly from the local stress.\\
Numerical computations suggest that optimal shapes are characterized by regions with substantially different mircostructure and sharp transitions between these regions. Thus an adaptive mesh strategy based on an a posteriori error control is evidently the appropriate strategy.
To this end, the error in the achieved cost has to be controlled.

{\it Review of related work.}
Shape optimization is a classical field in PDE constrained optimization and has been covered extensively in the literature, see \eg\
the textbooks \cite{Be95a,Al02}. For the link to homogenization theory we refer to \cite{Jikov1994,BrDe98,DoCi99,BuDa91}. 
Optimal microstructures where 
first derived by Hashin in 1962 via the concentric sphere construction for hydrostatic loads \cite{Ha62}. The sequential lamination construction
dates back to the 1980s \cite{Ta85,MuTa85,FrMu86,LuCh86,GiCh87,Av87} and was later used in a practical numerical scheme for topology optimization
in~\cite{AlBoFr97}. In all these cases proofs of optimality rely on the Hashin-Shtrikman bounds on the attainable sets of effective elastic properties \cite{HaSh63}.
Related to the homogenization approach is the so called free material optimization, see \eg~\cite{HaKoLe10}. Here the coefficients of the elasticity
are the degrees of freedom to be optimized. Possibly yielding materials without a physical equivalent, best approximating realizations are searched
in a post processing step, sometimes referred to as inverse homogenization. 
Applying the classical homogenization paradigm one asks for the optimal (geometric) microstructures and the effective material properties 
such that a macroscopic cost functional is optimized.
In~\cite{BaTo10,BaTo10a} a displacement approach on the microscale is investigated which consists of affine plus periodic functions on
generalized periodic domains. For the numerical optimization of the microscopic geometry a boundary tracking method was used. 
The incorporation of topological derivatives allowed for the nucleation of small holes and additional remeshing steps ensured mesh quality.
For specific geometric microstructures the effective stress response of
basic affine strains on the microscopic cell boundary are investigated in \cite{BaTo12} and used to describe the effective macroscopic behavior. 
The geometry of the locally periodic microstructures are then optimized.
The same methodology was picked up in \cite{CoGeRu14} while investigating microstructured materials with geometrically simple
perforations. The resulting shapes offered compliance values close to optimal ones generated by twofold sequential lamination. Moreover they
indicate the need to resolve sharp interfaces between regions with different types of local microstructure. 
To be able to both capture microscopic effects and macroscopic material behavior while maintaining feasible computational complexity in 
PDE simulation tailored numerical methods have to be proposed.
The heterogeneous multi-scale method (HMM) presented in \cite{EEnHu03,EEn03,EEn05,EMiZh05} depicts a very general
paradigm using independent macroscopic and microscopic schemes.
In \cite{Oh05,HeOh09} a posteriori error estimates for elliptic homogenization problems involving a fine scale diffusion term were derived.
Here the heterogeneous multi-scale method is reformulated as a direct finite element discretization of the two-scale homogenized
equation in variational form.
The obtained error indicator provides separate terms for errors in the macroscopic and microscopic discretization.

Concerning the a posteriori control  it is often not desirable to measure errors in classical function space norms but the error for a prescribed cost functionals.
From a practical point of view quantities of interest for a real composite work piece like \eg\ stresses in critical sections,
averaged surface tensions or (mollified) pointwise stresses can be considered as cost functionals and corresponding a
posteriori error estimates are derived \cite{PrOd99,OdVe00,OdVe00a,Ve04}. \\
The approach we pick up in this contribution is the dual weighted residual (DWR) method introduced in \cite{BeRa97} and applied to optimal control problems 
for instance in \cite{BeKaRa00}. Here residual type error estimates get weighted using duality methods.
In \cite{BeEsTr11} a special marking strategy was employed allowing to show quasi-optimality of the associated adaptive finite element
method while maintaining sound convergence. 
Control and state constraints were taken into account in \cite{BeVe09,VeWo08,LeMeVe13}.
Matrix valued $L^1$ optimal control in coefficients is studied in \cite{KoLe13}.

Shape optimization was addressed in \cite{KiVe13} using a tracking type cost functional, a Helmholtz state equation, and a graph representation of
the boundary. For this control function higher regularity could be shown and a priori estimates were derived.
An adaptive finite element method for shape optimization via boundary tracking, remeshing and perimeter penalization was presented in
\cite{MoNoPa10,MoNoPa12}. Here a DWR type of approach is used to assess the PDE error and the actual geometric error is estimated differently.
Furthermore the DWR method was used in \cite{KaDeGa14, Ka13} in the context of one shot methods for fuel ignition problems, the viscous Burgers
equation and aerodynamic shape optimization.
Optimal design in Navier Stokes flow is considered in \cite{BrLiUl12}. Here the DWR method for optimal control problems leads to separate terms for the
primal, the adjoint and the control residual. 
Similarly \cite{Wo10} addresses a problem in free material optimization based on the variable thickness sheet model.

The paper is organized as follows. In section~\ref{sec:basics} we will describe in condensed form the necessary background on linearized elasticity, shape optimization, the sequential laminates construction in shape optimization, and its numerical discretization.
A posteriori error estimates based on the dual weighted residual approach are derived in section~\ref{sec:dwr}.
In section~\ref{sec:apriori} we discuss the estimation of weighting term using a priori regularity
while section~\ref{sec:localapprox} focuses on their numerical approximation. Details concerning the implementation
are given in section~\ref{sec:impl} and numerical results are finally presented in section~\ref{sec:results}.

\section{Elastic shape optimization and microscopic lamination}\label{sec:basics}
We consider an elastic object given as a simple connected domain $\body \subseteq \workdom \subset \R^2$
where $\workdom$ is a circumjacent working domain and we at first suppose that $\workdom$ is filled with elastic material (hard and soft)
and we aim at the optimization of the distribution of the two phases. 
A Dirichlet boundary part $\Gamma_D \subset \partial D \cap \partial \body$ and a Neumann boundary 
$\Gamma_N \subset \partial D \cap \partial \body$ subject to a sufficiently regular surface load $g$ are both supposed to be fixed and not subject to optimization. 
The remaining boundary $\partial \Omega \setminus (\Gamma_D \cup \Gamma_N)$ can be varied. 
The applied forcing causes stresses inside the material which maintain an inner equilibrium. It can be described
by the system of partial differential equations  of linearized elasticity. Let $\Sigma$ denote the set of 
feasible stresses in accordance to the PDE, i.e.
\begin{equation}\label{eq:feasible}
 \Sigma\! :=\! \left\lbrace \stress\!:\! \workdom \!\rightarrow\! \R^{2\times 2}_{\text{sym}} \; \big| \;\!
 \div \{ \stress \} \!= \!0 \text{ in }\workdom, \, \stress n \!=\! g \text{ on }\Gamma_{\!N}, \, \strain \!=\! 0 \text{ on }\Gamma_{\!D},
 \, \stress \!= \!\etensor[q]  \eps{\strain} \right\rbrace
\end{equation}
Here $\strain: \workdom \rightarrow \R^2$ denotes the displacement, $n$ the outward pointing normal on $\Gamma_N$,
$\eps{\strain}  = \frac12 \left( \D\strain + \D\strain^\top \right)$ the symmetrized strain tensor, and
$\etensor[q]: \workdom \rightarrow \R^{2^4}$ the elasticity tensor. The elasticity  tensor is a forth order tensor with the symmetry properties
$(\etensor[q])_{ijkl} = (\etensor[q])_{jikl} = (\etensor[q])_{ijlk} = (\etensor[q])_{klij}$ and the ellipticity condition
$\etensor[q] \, \xi : \xi \geq c \, \xi \cdot \xi$. Usually it can be characterized locally by a small parameter vector $q$.
We consider here material which behaves  isotropic  and is described by  Lam\'e parameters $\lambda$ and $\mu$ and the strain stress relation 
$\stress = 2 \mu \eps{\strain} + \lambda \tr\eps{\strain} \1$.\\
For later purposes we take into account the weak displacement formulation of \eqref{eq:feasible}, \ie\
$\strain \in H^1_{\Gamma_D}$, the space of $L^2$ vector valued functions with weak derivatives in $L^2(\Omega)$ 
and vanishing trace on $\Gamma_D$, solves the variational problem
\begin{equation}\label{eq:weak}
 a(q;\strain,\varphi) = l(\varphi) \;\;\forall\,\varphi \in H^1_{\Gamma_D} 
\end{equation}
with the quadratic form $a(q;\strain,\varphi) = \int_\workdom \etensor[q] \, \eps{\strain} : \eps{\varphi} \dx$ and the linear form
$l(\varphi)     = \intga \etensor[q] \, \eps{\strain} n \cdot \varphi \dSx = \intga g \cdot \varphi \dSx \,$.

\begin{figure}[!ht]
\hfill
\begin{minipage}{.22\linewidth}
\setlength{\unitlength}{\linewidth}
\includegraphics[width=\linewidth]{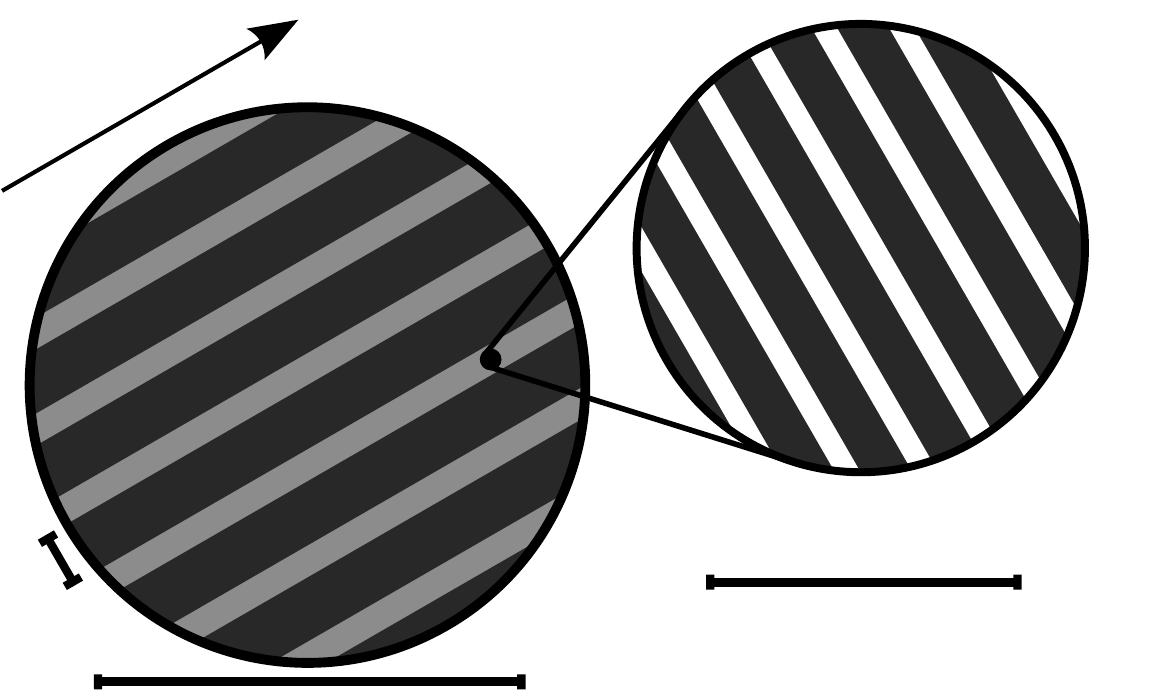}
\begin{picture}(0,0)
\put(.08,.68){$\alpha$}
\put(-.05,.19){$m$}
\put(.71,.13){$\epsilon_1$}
\put(.22,.04){$\epsilon_2$}
\end{picture}
\end{minipage}
\hfill
\begin{minipage}{.25\linewidth}
\includegraphics[width=\linewidth]{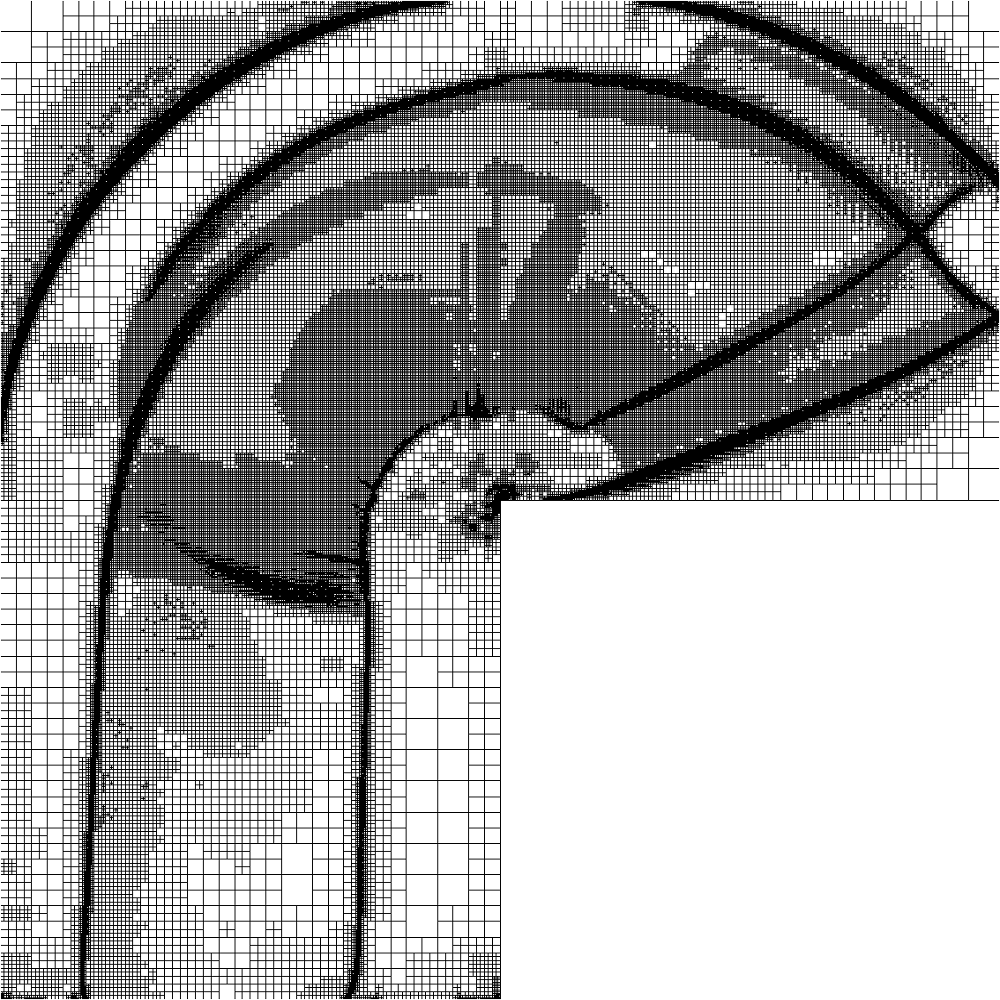}
\end{minipage}
\hfill
\begin{minipage}{.25\linewidth}
\includegraphics[width=\linewidth]{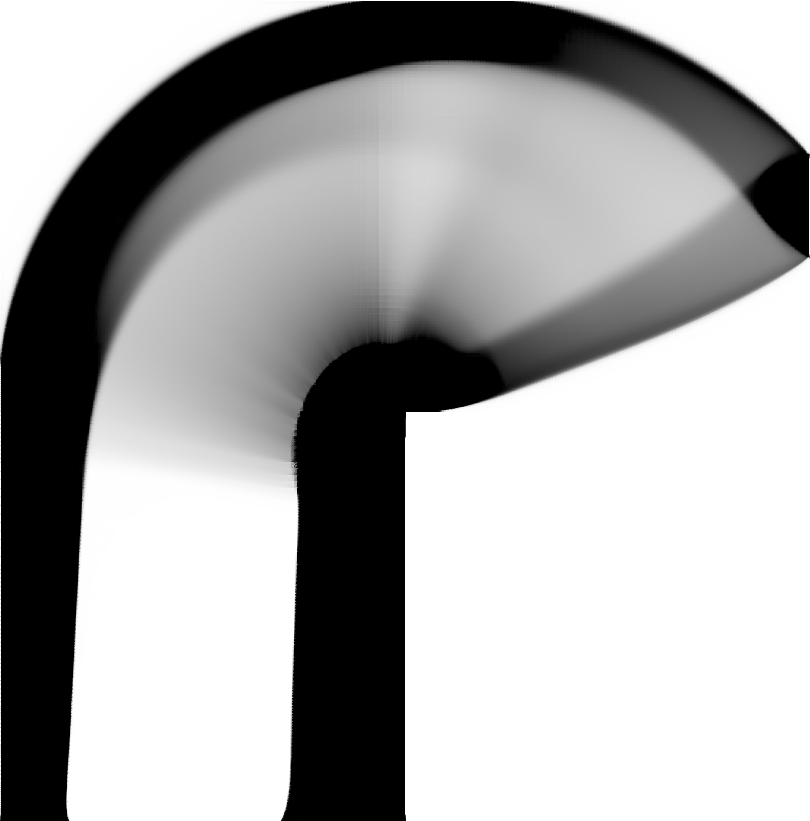}
\end{minipage}
\hfill
\begin{minipage}{.25\linewidth}
\includegraphics[width=\linewidth]{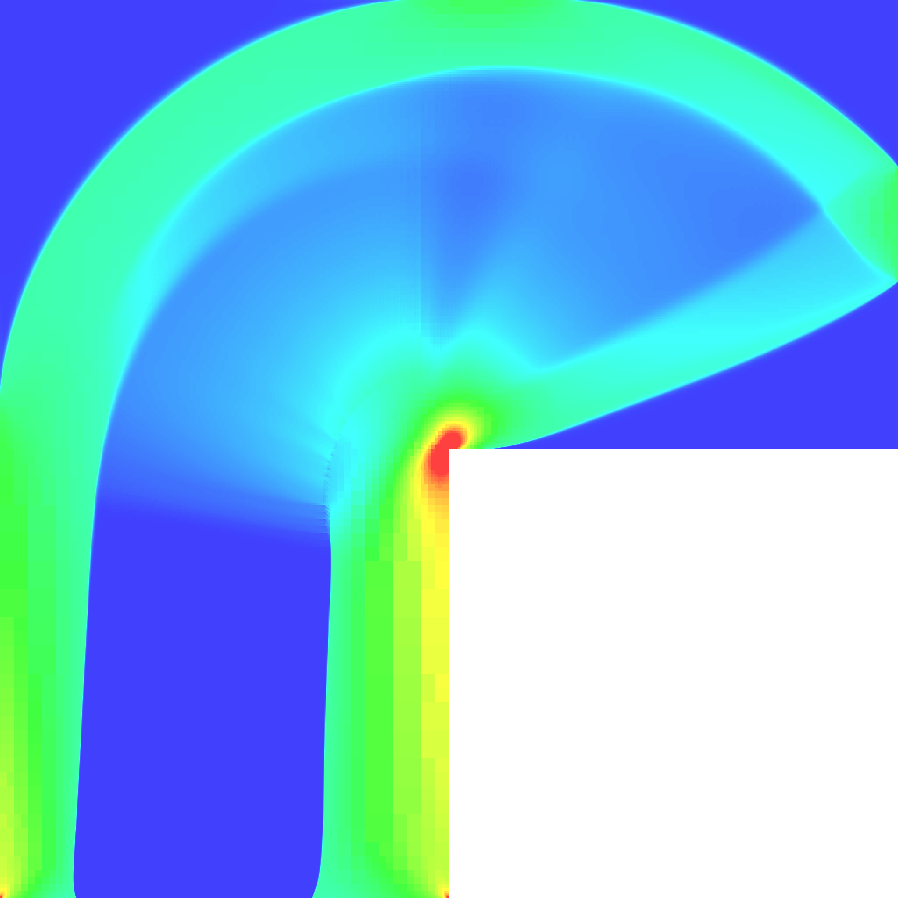}
\end{minipage}
\hfill\mbox{}
\caption[]{Sketch of twofold sequential lamination.
           Numerically optimized result for an L-shaped domain fixed at the bottom with downwards pointing load at the center region on the right hand side.
           An adaptively refined grid is shown along with the density lamination parameter $\theta$ plotted and the von Mises stress using the color coding \raisebox{1pt}{\includegraphics[height=4pt]{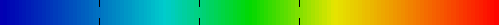}}.}
\label{fig:laminatesSketch}
\end{figure}
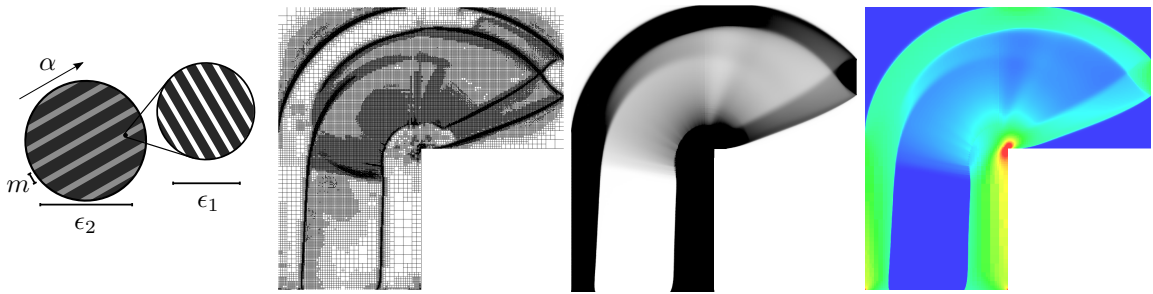

A classical shape optimization problem now amounts to the task of distributing a hard material whose elastic behavior is
described by an elasticity tensor $A$. To this end, we consider a characteristic function $\chi \in L^\infty(\workdom,\lbrace 0,1 \rbrace)$ of the hard phase
with $\int_\workdom \chi \dx = \Theta$ reflecting the volume constraint. The elasticity tensor on $\workdom$ is then
defined as $\etensor[q] = \chi A[q] + (1-\chi) B[q]$ where $B[q]$ is an elasticity tensor for a soft material filling 
the remaining space.
To assess the performance of a given shape for a prescribed loading a suitable cost functional has to be taken into account.
We restrict ourselves to the compliance cost as a global measure of rigidity. For a displacement $\strain= \strain[\chi]$ being the unique solution
of \eqref{eq:weak} for given $\chi$ it is given by
\begin{equation} \label{eq:compliance}
 J[\strain[\chi] ,\chi ] := \intga g \cdot \strain[\chi] \dSx = l(\strain[\chi] ) = a(q;\strain[\chi] ,\strain[\chi] ) \,.
\end{equation}
The cost functional $ J[\strain[\chi] ,\chi ]$ should thereby be minimized with respect to $\chi$ subject to a constraint on the volume of the hard material phase.
 Ultimately one is interested in the degenerate case $B=0$, where there is void outside of 
$\Omega=\{x\in \workdom \;|\; \chi(x) = 1\}$.
This minimizing problem is ill-posed and on a minimizing sequence the onset of microstructures can be observed. 
Thus one investigates a relaxed formulation of this shape optimization problem. We refer to
\cite{AlBoFr97,Al02} for a concise summary and references therein for further details. Relaxation leads to the following reformulation:
\begin{equation} \label{eq:relaxed}
\min\limits_{\stress \in \Sigma} \, \int_\workdom \min\limits_{0 \leq \theta \leq 1} \,
\min\limits_{\efftensor[q] \in G_\theta^B} {\efftensor[q]}^{-1} \stress : \stress + l \, \theta \dx.  
\end{equation}
Reading from left to right, for fixed stresses
$\sigma$ on $\workdom$ solving \eqref{eq:feasible}, a fixed point $x \in D$ and a fixed local density $\theta(x)$, an
elasticity tensor $\efftensor[q]$ minimizing the elastic energy density is to be found 
among all effective homogenized elasticity tensors resulting from a mixture of materials $A$ and $B$.
Additionally the Lagrangian multiplier $l$ takes into account the amount of total material used. It can be
shown that it is a decreasing function of $l$ which facilitates an adaptation to enforce a fixed prescribed
volume at any time.
For the set $G_\theta^B$ of homogenized tensors, also referred to as G-closure, no algebraic characterization is known. However there exist
sharp bounds on the elastic energy density, called Hashin-Shtrikman bounds. Theses bounds are attained within the subclass
of sequentially laminated materials. The construction starts on a scale $\epsilon_1$ by layering rigid material $A$ and
soft material $B$ with a certain ratio and along a certain direction. By means of homogenization effective material properties
can be computed and used in the next iteration in place of material $B$. Mechanically this takes place on a substantially larger scale $\epsilon_2 \gg \epsilon_1$.
In two dimensions and for the compliance cost it is known that the construction on these two iterations is sufficient, cf. \cite{KoLi88,AlKo93a}.
It is moreover possible to pass from the weak material $B$ to void \cite{AlKo93,AlBoFr97}. \\
The sequential laminates microstructure is locally characterized by three parameters (cf.~Fig.~\ref{fig:laminatesSketch}):
the inclination of the main lamination direction $\alpha(x)$, the ratio of material spent in the second lamination stage $m(x)$,
and the overall local material density $\theta(x)$. The second lamination direction is always orthogonal to the first one and the second
ratio parameter is $1-m(x)$.

All three parameters directly depend on the eigenvalues $\lambda_1$, $\lambda_2$
and the eigenvectors $\tau^{(1)}$, $\tau^{(2)}$ of the local stress tensor $\stress$ and can be computed explicitly via the mapping
$x \mapsto q[\stress](x) = \left( \alpha[\stress](x), m[\stress](x),\theta[\stress](x) \right)^\top$, see \cite{AlBoFr97}:
\begin{equation}\label{eq:explicitparameters}
\begin{aligned}
 \alpha[\stress] &= \arctan\left( \frac{\tau^{(1)}_y(\stress)}{\tau^{(1)}_x(\stress)} \right) \,, \qquad
 m[\stress]      = \frac{|\lambda_2(\stress)|}{|\lambda_1(\stress)|+|\lambda_2(\stress)|} \,, \\
 \theta[\stress] &= \min\left\{ 1, \sqrt{ \frac{2 \mu + \lambda}{4 \mu (\mu+\lambda) \, l} } ( |\lambda_1(\stress)|+|\lambda_2(\stress)| ) \right\} \,.
\end{aligned}
\end{equation}
Here $\lambda$ and $\mu$ denote the Lam\'e parameters of the rigid isotropic material $A$ and $\tau^{(i)}_x$, $\tau^{(i)}_y$ the first and
second component of the two dimensional vector $\tau^{(i)}$.
Given these parameters the effective homogenized elasticity tensor $\efftensor[q](x)$ can likewise be computed explicitly:

\begin{equation}\label{eq:effectivetensor}
\begin{aligned}
C^\ast_{mnop}[q] &= R[\alpha] \, \bar{C}[m,\theta] := Q_{mi}[\alpha] \, Q_{nj}[\alpha] \, Q_{ok}[\alpha] \, Q_{pl}[\alpha] \; \bar{C}_{ijkl}(m,\theta) \,, \\[6pt]
\bar{C}_{1111}[m,\theta] &= \frac{ 4 \kappa \mu (\kappa+\mu) \theta (1-\theta (1-m))(1-m) } {4 \kappa \mu m (1-m) \theta^2 + (\kappa+\mu)^2 (1-\theta)} \,, \\
\bar{C}_{2222}[m,\theta] &= \frac{ 4 \kappa \mu (\kappa+\mu) \theta (1-\theta m)m }         {4 \kappa \mu m (1-m) \theta^2 + (\kappa+\mu)^2 (1-\theta)} \,, \\
\bar{C}_{1122}[m,\theta] &= \frac{ 4 \kappa \mu \lambda \theta^2 m (1-m) }                  {4 \kappa \mu m (1-m) \theta^2 + (\kappa+\mu)^2 (1-\theta)} \,.
\end{aligned}
\end{equation}
Thereby $R$ denotes the linear mapping which transforms the tensor $\bar{C}$ given in reference configuration into the appropriate coordinate frame
spanned by the eigenvectors of the stress tensor. In its definition $Q$ are $2 \times 2$ rotation matrices and the Einstein
summation convention is used. The bulk modulus is defined as $\kappa  = \lambda + \mu$. The tensor $\bar{C}$ is complemented using the symmetry relations and filling the remaining entries with  $0$. This yields a singular elasticity tensor in the sense that a strain with
just off diagonal entries, corresponding to shearing deformations, will lead to a vanishing stress. \\
In the following we will assume that the density $\theta$ is always bounded away from $0$ by a small threshold.
Thus we deal with the hard\/soft configuration. For the ease of notation
we extend the Neumann boundary to $\Gamma_N = \partial D \setminus \Gamma_D$ and prescribe
homogeneous boundary conditions $g = 0$ on the added parts.

For the discretization we use a finite element ansatz and assume that we are given a 
triangulation $\mathcal{T}$ of $\workdom$ with elements $T \in \mathcal{T}$ and denote by $h$ the piecewise constant mesh size function.
 Let $\Vh$ be the space of piecewise bi-quadratic and continuous vector valued functions on which  we ask for a solution $\strain_h$ of the discrete
weak problem $a(q_h;\strain_h,\varphi_h) = l(\varphi_h)$ for all $\varphi_h \in \Vh$, or more explicitly
\begin{equation}\label{eq:discrete}
  \begin{aligned}
 & \int_\workdom \efftensor[q_h] \, \eps{\strain_h} : \eps{\varphi_h} \dx
                 = \intga g \cdot \varphi_h \dSx  \;\;\forall\,\varphi_h \in \Vh \,.
  \end{aligned}
\end{equation}
Associated to $\strain_h$ are stresses $\stress_h := \stress[u_h](x) = \efftensor[q_h] \, \eps{u_h}$
for which lamination parameters are computed via formula \eqref{eq:explicitparameters}.
For brevity of notation we write $q_h = (\alpha_h,m_h,\theta_h) := q[\stress_h]$.
Let us emphasize that the parameter function varies from point to point.
More details on the numerical realization will be given in section~\ref{sec:impl}.

\section{A posteriori error estimates based on the dual weighted residual approach} \label{sec:dwr}
In this section we will derive estimates of the error in the compliance objective \eqref{eq:compliance} when using optimal 
sequential laminated microstructures in shape optimization. Let $\strain$ be the solution of the continuous problem \eqref{eq:weak}
involving the elasticity tensor $\efftensor[q]$ with optimal lamination parameters $q$ defined by \eqref{eq:explicitparameters}
and $\strain_h$ be a discrete solution satisfying \eqref{eq:discrete}.

We are interested in the a posteriori control of the numerical approximation error in the cost functional.
To this end we pick up the dual weighted residual approach  \cite{BeKaRa00}. In that respect we will treat 
the rotation parameter $\alpha$ differently from the other two parameters $m$ and $\theta$, which will from now on be denoted by
$\vartheta := (m,\theta)$. 

Let us begin with a useful observation for the elastic energy density, namely that 
it is invariant with respect to the rotation of the microstructure. Intuitively this follows from
the fact that the lamination construction yields an orthotropic material which is always aligned with the main stress directions.

\begin{lemma}[Invariance of the elastic energy density \wrt\  rotations] \label{lemma:invariance}
Let $x \in \workdom$ be fixed and $\efftensor[\alpha,\vartheta]$ be an effective elasticity tensor corresponding to an optimal
rank 2 sequential laminate at $x$. Let $\strain[\alpha,\vartheta]$ be the solution of \eqref{eq:weak}.
Then the local elastic energy density is invariant \wrt\  the rotation parameter $\alpha$,  \ie
$$\tdiff{\alpha}  \left( \efftensor[\alpha,\vartheta] \, \eps{\strain[\alpha,\vartheta]} : \eps{\strain[\alpha,\vartheta]} \right) (x) = 0\,.$$
\end{lemma}

\textbf{Proof.}
With $\sigma[\alpha,\vartheta] := \efftensor[\alpha,\vartheta] \, \eps{\strain}$ the elastic energy density can be rewritten as
\[
 \efftensor[\alpha,\vartheta] \, \eps{\strain[\alpha,\vartheta]} : \eps{\strain[\alpha,\vartheta]} = {\efftensor}^{-1}[\alpha,\vartheta] \, \sigma[\alpha,\vartheta] : \sigma[\alpha,\vartheta] \,.
\]
Here $\alpha$ and $\vartheta$ enter the effective homogenized tensor $\efftensor[\alpha,\vartheta]$ by means of \eqref{eq:effectivetensor}.
Since $\efftensor[\alpha,\vartheta]$ represents an optimal material the elastic energy density attains the lower Hashin-Shtrikman bound,
cf.~\cite[Theorem~2.3.35]{Al02}, and we obtain 
\[
 {\efftensor}^{-1}[\alpha,\vartheta] \, \sigma[\alpha,\vartheta] : \sigma[\alpha,\vartheta]
 = A^{-1} \sigma[\alpha,\vartheta] : \sigma[\alpha,\vartheta]+ \frac{(\kappa+\mu)\theta}{4\kappa \mu (1-\theta)} ( |\lambda_1| + |\lambda_2| )^2 \,.
\]
Here the eigenvalues $\lambda_1$, $\lambda_2$ of $\sigma$ do not depend on the rotation.
In fact, the rotation parameter $\alpha$ is derived from the orientation of the eigenvectors.
Thus $\alpha$ only appears via $\sigma$ in the first term of the right hand side.
This term also represents an elastic energy density but this time involving the isotropic constituent $A$.
Because of the isotropy of $A$ altering the alignment of $\sigma$, \eg\ by rotating to reference configuration, leaves this term unchanged. 
Altogether this means that the right hand side does not depend on the rotation parameter. This proofs the claim.
\qed \\

We now present the main result.
\begin{theorem}[Weighted a posteriori error estimate] \label{theorem:estimates}
For a given continuous solution $\strain$ to \eqref{eq:weak} with a microstructured material given by the optimal lamination parameter function $q$ and a numerical
approximation $u_h$ the following error estimate holds for the compliance objective \eqref{eq:compliance}:
\[
 |J[q;u]-J[q_h;u_h]| \leq \sum_T \eta_T(u_h,q_h) + \mathcal{R} \,.
\]
Here $\mathcal{R}$ is a higher order remainder term involving higher order derivatives and the local error indicators $\eta_T$ of first order are defined by
\[
  \eta_T(u_h,q_h) :=  \rho_T^{(u)} \omega_T^{(u)} + \rho_{\partial T}^{(u)} \omega_{\partial T}^{(u)}
                   + \fr12 \rho_T^{(m)} \omega_T^{(m)} + \fr12 \rho_T^{(\theta)} \omega_T^{(\theta)} 
  \quad \text{with}
\]
\begin{equation*}
\begin{aligned}
  \rho_T^{(u)}              &= \left\Vert -\div \left\{ \stress_h \right\} \right\Vert_{0,2,T} \,, &
  \omega_T^{(u)}            &= \left\Vert u-\tilde{u}_h \right\Vert_{0,2,T} \,,  \\
  \rho_{\partial T}^{(u)}   &= 
  \left\lbrace 
  \begin{array}{l l}
    \Vert \textstyle\frac{1}{2} \left[ \stress_h n \right]\Vert_{0,2,\partial T}; & \partial T \cap \partial\workdom = \emptyset \\
    \Vert \stress_h n - g\Vert_{0,2,\partial T};                                 & \partial T \subset \Gamma_N
  \end{array}
  \right. \,,   &
  \omega_{\partial T}^{(u)} &= \left\Vert u-\tilde{u}_h \right\Vert_{0,2,\partial T} \,,  \\
  \rho_T^{(m)}         &= \left\Vert R[\alpha_h] \, \bar{C}_{,m}\left[m_h,\theta_h\right] \; \eps{\strain_h[\alpha_h]} : \eps{\strain_h[\alpha_h]}\right\Vert_{0,1,T} \,,  &
  \omega_T^{(m)}       &= \left\Vert m-\tilde{m}_h \right\Vert_{0,\infty,T} \,,  \\
  \rho_T^{(\theta)}    &= \left\Vert R[\alpha_h] \, \bar{C}_{,\theta}\left[m_h,\theta_h\right] \; \eps{\strain_h[\alpha_h]} : \eps{\strain_h[\alpha_h]}\right\Vert_{0,1,T} \,,  &
  \omega_T^{(\theta)}  &= \left\Vert \theta-\tilde{\theta}_h \right\Vert_{0,\infty,T} \,,
\end{aligned}
\end{equation*}
where  $\tilde{u}_h \in \Vh$ and $\tilde{m}_h$, $\tilde{\theta}_h$ are chosen arbitrarily.
\end{theorem}

\textbf{Proof.} 
The first step is to use the Lagrangian $\Lag(q;u,p) := J[q;u] + a(q;u,p) - l(p)$  for decoupled $q$, $u$, and $p$ and take into account 
the PDE constraint, namely that the displacement $\strain$ has to be a solution of the linearized elasticity system \eqref{eq:weak}. 
The first order optimality conditions are
\begin{align}
 \rho(q;u)(v)        := \Lag_{,p}(q;u,p)(v) &= a(q;u,v) - l(v) \,,  \label{primal}\\
 \rho^\ast(q;u,p)(v) := \Lag_{,u}(q;u,p)(v) &= J_{,u}[q;u](v) + a(q;v,p) \,,  \label{dual} \\
 \rho^q(q;u,p)(v)    := \Lag_{,q}(q;u,p)(v) &= J_{,q}[q;u](v) + a_{,q}(q;u,p)(v) - l_{,q}(p)(v) \,,  \label{control}
\end{align}
where $v$ is a universal variable for an arbitrary test function from the corresponding space.
Equation \eqref{primal} is just the definition of the elastic solution $u$.
Equation \eqref{dual} defines the dual solution, which in case of compliance is  $p=-u$.
As we assumed $u$ and $u_h$ to be solutions of the continuous and the discrete problem, respectively,
$(q,u,-u)$ and $(q_h,u_h,-u_h)$ are stationary points of the Lagrangian. Thus we have
\[
 \Lag(q;u,p)-\Lag(q_h;u_h,p_h) = J[q;u]-J[q_h;u_h] \,.
\]
Estimating the error can therefore be done by considering the difference in the Lagrangian.\\

Next we consider a suitable first order expansion of the Lagrangian. 
Here we will take special care of the rotation parameter $\alpha$ to be able
to use Lemma~\ref{lemma:invariance} above. 
The Lagrangian thus only depends on $u$ and the controls $q=(\alpha,\vartheta)$ because the dual solution coincides with the negative primal solution.
Hence, in the following we will skip the dual solution as a parameter. 
We now study the dependence  of the elastic solution $\strain[\alpha]$ on the rotation parameter $\alpha$.\\
As a shortcut notation we write $e_u(s) = \strain[\alpha]-\strain_h[\alpha_h + s e_\alpha]$, $e_\alpha = \alpha -\alpha_h$, and $e_{\vartheta} = \vartheta-\vartheta_h$
for the difference between the exact and the discrete displacement and controls. Here $u_h$ is given for the continuous spatially varying parameter function
$\alpha_h + s e_\alpha$ as the solution to the corresponding problem \eqref{eq:discrete}.
We now express the error in the Lagrangian as a 1D integral and observing $\left. \strain_h[\alpha_h+s e_\alpha] + s e_u(s) \right|_{s=1} =  \strain[\alpha]$ 
and $\left. \strain_h[\alpha_h+s e_\alpha] + s e_u(s) \right|_{s=0}$ $= \strain_h[\alpha_h]$ we obtain
\begin{eqnarray*} 
  e_\Lag &:=&  \Lag(\alpha,\vartheta;\strain[\alpha])-\Lag(\alpha_h,\vartheta_h;\strain_h[\alpha_h])\\
  &=& \int_0^1 \tdiff{s} \Lag (\alpha_h + s e_\alpha, \vartheta_h + s e_{\vartheta},  \strain_h[\alpha_h+s e_\alpha] + s e_u(s) ) \,\mathrm{d}s\,.
\end{eqnarray*}
Now, we use the trapezoidal rule $\int_0^1 f(s) \mathrm{d}s  =  \frac12 ( f(0) + f(1) ) + \frac12 \int_0^1 f''(s) s(1-s)\, \mathrm{d}s$ with
$f(s) = \tdiff{s} \Lag(\alpha_h + s e_\alpha, \vartheta_h + s e_{\vartheta}, \strain_h[\alpha_h+s e_\alpha] + s e_u(s) )$.
The term $f(1)$ vanishes because $(\alpha,\vartheta,\strain[\alpha])$ is assumed to be a stationary point and
$e=(e_\alpha,e_{\vartheta},\strain[\alpha] - \strain_h[\alpha])^T$ is a feasible test direction, i.e. $\nabla \Lag(\alpha,\vartheta,\strain[\alpha])(e) = 0$.
Thus, we obtain 

$$
e_\Lag = \frac12  \nabla \Lag (\alpha_h, \vartheta_h, \strain_h[\alpha_h])
                  \left( \begin{array}{c} e_\alpha \\ e_{\vartheta} \\ {\strain_h}_{,\alpha}[\alpha_h](e_\alpha) + e_u(0) \end{array} \right) + \mathcal{R}\,
$$
where $\mathcal{R} := \frac12 \int_0^1 \frac{\mathsf{d}^3}{\mathsf{d} s^3} \Lag
  ( \alpha_h + s e_\alpha, \vartheta_h + s e_{\vartheta}, \strain_h[\alpha_h+s e_\alpha] + s e_u(s)) \, s \, (1-s) \mathrm{d}s\,$.
Next, we expand the lower order term and achieve
\begin{eqnarray*}
&& \frac12 \nabla \Lag (\alpha_h, \vartheta_h, \strain_h[\alpha_h])
                  \left( \begin{array}{c} e_\alpha \\ e_{\vartheta} \\ {\strain_h}_{,\alpha}[\alpha_h](e_\alpha) + e_u(0) \end{array} \right) \notag \\
    && = \frac12 \Lag_{,u}         (\alpha_h,\vartheta_h;\strain_h[\alpha_h]) ( \strain[\alpha] - \strain_h[\alpha_h] )
      + \frac12 \Lag_{,\vartheta} (\alpha_h,\vartheta_h;\strain_h[\alpha_h]) ( \vartheta-\vartheta_h ) \notag  \\
    && \quad + \frac12 \Lag_{,\alpha}    (\alpha_h,\vartheta_h;\strain_h[\alpha_h]) ( \alpha -\alpha_h )
      + \frac12 \Lag_{,u}         (\alpha_h,\vartheta_h;\strain_h[\alpha_h]) ( {\strain_h}_{,\alpha}[\alpha_h](\alpha -\alpha_h))   \notag  \\
    && = \frac12 \Lag_{,u}         (\alpha_h,\vartheta_h;\strain_h[\alpha_h]) ( \strain[\alpha] - \strain_h[\alpha_h] )
      + \frac12 \Lag_{,\vartheta} (\alpha_h,\vartheta_h;\strain_h[\alpha_h]) ( \vartheta-\vartheta_h ) \notag  \\ 
    && \quad  + \frac12 \tdiff{\alpha} \Lag (\alpha,\vartheta_h;\strain_h[\alpha])\big|_{\alpha=\alpha_h} ( \alpha -\alpha_h ) \,. \notag  
\end{eqnarray*}
By definition we have
\[
    \tdiff{\alpha} \Lag (\alpha,\vartheta_h;\strain_h[\alpha])\big|_{\alpha=\alpha_h} 
  = \left[\tdiff{\alpha} \int_\workdom \efftensor\left[\alpha,\vartheta_h \right] \eps{\strain_h[\alpha] } : \eps{\strain_h[\alpha]} \dx\right]_{\alpha=\alpha_h}
\]
and thus the term vanishes by virtue of Lemma~\ref{lemma:invariance}.
Finally we consider arbitrary discrete test functions $\tilde{u}_h \in \Vh$ and $\tilde{\vartheta}_h$
resulting from an arbitrary stress by means of formula \eqref{eq:explicitparameters}
and use the fact that $(\alpha_h,\vartheta_h,\strain_h[\alpha_h])$ is a discrete stationary point of the Lagrangian and thus $\nabla \Lag (\alpha_h, \vartheta_h, \strain_h[\alpha_h])$ vanishes in directions $(0, \tilde{\vartheta}_h -\vartheta_h, \tilde{\strain}_h[\alpha_h] - \strain_h[\alpha_h] )$. This gives
\begin{eqnarray*}
e_\Lag \!\!
&=& \!\! \frac12 \Lag_{,u}         (\alpha_h,\vartheta_h;\strain_h[\alpha_h]) ( \strain[\alpha] - \tilde{\strain}_h[\alpha_h] )
+ \frac12 \Lag_{,\vartheta} (\alpha_h,\vartheta_h;\strain_h[\alpha_h]) ( \vartheta - \tilde{\vartheta}_h ) + \mathcal{R} \,.\notag
\end{eqnarray*}
Furthermore, we consider the weak formulation \eqref{eq:weak} and the compliance cost functional
\eqref{eq:compliance}. For the \emph{primal residual} we get:
\begin{eqnarray*}
   && \Lag_{,u} (\alpha_h,\vartheta_h;\strain_h[\alpha_h]) ( \strain[\alpha] - \tilde{\strain}_h[\alpha_h] ) \\
 &&=  -2 \int_\workdom \efftensor\left[\alpha_h,\vartheta_h \right] \eps{\strain_h[\alpha_h] } : \eps{ \strain[\alpha] - \tilde{\strain}_h[\alpha_h]} \dx + 2 \intga g \cdot (\strain[\alpha] - \tilde{\strain}_h[\alpha_h]) \\
 &&= 2 \, \Big( \sum_j \int_{T_j} \div \left\{ \efftensor\left[\alpha_h,\vartheta_h \right] \eps{\strain_h[\alpha_h]} \right\} \cdot (\strain[\alpha] - \tilde{\strain}_h[\alpha_h]) \dx \\
   && \qquad \qquad - \int_{\partial T_j} \efftensor\left[\alpha_h,\vartheta_h \right] \eps{\strain_h[\alpha_h]} n \cdot (\strain[\alpha] - \tilde{\strain}_h[\alpha_h]) \dSx \\
   && \qquad \qquad + \int_{\partial T_j \cap \Gamma_N} g \cdot (\strain[\alpha] - \tilde{\strain}_h[\alpha_h]) \dSx \Big)
\end{eqnarray*}
This leads to the postulated residual terms with 
\[
\begin{aligned}
 \rho_T^{(u)}              &= \left\Vert \div \left\{ \stress[\strain_h[\alpha_h]] \right\} \right\Vert_{0,2,T}, &  \!\!
 \omega_T^{(u)}            &= \left\Vert \strain[\alpha] - \tilde{\strain}_h[\alpha_h] \right\Vert_{0,2,T}, \\
 \rho_{\partial T}^{(u)}   &=
 \left\lbrace
  \begin{array}{l l}
   \Vert \textstyle\frac{1}{2} \left[ \stress[\strain_h[\alpha_h]] \cdot n \right]\Vert_{0,2,\partial T}\! & \!\partial T_j \cap \partial\workdom = \emptyset \\
   \Vert \stress[\strain_h[\alpha_h]] \cdot n - g\Vert_{0,2,\partial T}                                 \! & \!\partial T_j \subset \Gamma_N
  \end{array},
 \right. &  \!\!
 \omega_{\partial T}^{(u)} &= \left\Vert \strain[\alpha] - \tilde{\strain}_h[\alpha_h] \right\Vert_{0,2,\partial T }.
\end{aligned}
\]
For the \emph{control residual} we get using \eqref{eq:effectivetensor}
\begin{eqnarray*}
   && \Lag_{,\vartheta} (\alpha_h,\vartheta_h;\strain_h[\alpha_h]) ( \vartheta - \tilde{\vartheta}_h )  \\
 &&=  \sum_j \int_{T_j} C_{,\vartheta}\left[\alpha_h,\vartheta_h \right] ( \vartheta - \tilde{\vartheta}_h ) \eps{\strain_h[\alpha_h]} : \eps{\strain_h[\alpha_h]} \dx \\
 &&=  \sum_j \int_{T_j} R[\alpha_h] \, \bar{C}_{,m}\left[m_h,\theta_h\right] (m-\tilde{m}_h)                \; \eps{\strain_h[\alpha_h]} : \eps{\strain_h[\alpha_h]} \\
 && \qquad \qquad  +                 R[\alpha_h] \, \bar{C}_{,\theta}\left[m_h,\theta_h\right] (\theta-\tilde{\theta}_h) \; \eps{\strain_h[\alpha_h]} : \eps{\strain_h[\alpha_h]}  \dx
\end{eqnarray*}
Concerning the last expression we know that the controls $m_h$, $\theta_h$ being the laminate parameters are bounded.
Differentiation of the entries of the effective tensor \eqref{eq:effectivetensor} \wrt\ $m$ and $\theta$ again yields bounded expressions
within the attainable range of the parameters. The tensors $\bar{C}_{,m}$ and $\bar{C}_{,\theta}$ are therefore bounded and the
whole integrand is bounded in $L_1$.
We finally obtain the desired residual estimate with 
\[
\begin{aligned}
  \rho_T^{(m)}         &= \left\Vert R[\alpha_h] \, \bar{C}_{,m}\left[m_h,\theta_h\right] \; \eps{\strain_h[\alpha_h]} : \eps{\strain_h[\alpha_h]}\right\Vert_{0,1,T} \,, &
  \omega_T^{(m)}       &= \left\Vert m-\tilde{m}_h \right\Vert_{0,\infty,T} \,, \\
  \rho_T^{(\theta)}    &= \left\Vert R[\alpha_h] \, \bar{C}_{,\theta}\left[m_h,\theta_h\right] \; \eps{\strain_h[\alpha_h]} : \eps{\strain_h[\alpha_h]}\right\Vert_{0,1,T} \,, &
  \omega_T^{(\theta)}  &= \left\Vert \theta-\tilde{\theta}_h \right\Vert_{0,\infty,T} \,.
\end{aligned}
\]
This proofs the claim.
\qed

\section{Estimation of weights using a priori regularity} \label{sec:apriori}
For the error estimates of Theorem~\ref{theorem:estimates} the residual terms $\rho_T$ can be computed straightforwardly 
from the discrete solution $u_h$. The weights $\omega_T$, however, depend on the exact (unknown) solution $u$.
Using a priori known regularity those terms can be estimated. For the primal solution classical interpolation estimates exist.
For the laminate parameters it is possible to estimate the error by exploiting the explicit
relation to the primal solution \eqref{eq:explicitparameters}.
The following estimates especially ensure robustness of the error estimates for a sufficiently smooth continuous solution $u$:
\begin{corollary}[A priori estimates of the weights]
Under the assumption that $u \in W^{2,\infty}$ the following estimates for the weighting terms hold:
\begin{equation}
 \begin{aligned}
   \omega_T^{(u)}            &= \left\Vert u-\tilde{u}_h \right\Vert_{0,2,T}                && \lesssim h(T)            \; |u|_{1,2,\omega_T} \,, \\
   \omega_{\partial T}^{(u)} &= \left\Vert u-\tilde{u}_h \right\Vert_{0,2,\partial T}       && \lesssim h(T)^{\frac12}  \; |u|_{1,2,\omega_{\partial T}} \,, \\
   \omega_T^{(m)}            &= \left\Vert m-\tilde{m}_h \right\Vert_{0,\infty,T}           && \lesssim h(T) |\ln h(T)| \; |u|_{2,\infty,T} \,, \\
   \omega_T^{(\theta)}       &= \left\Vert \theta-\tilde{\theta}_h \right\Vert_{0,\infty,T} && \lesssim h(T) |\ln h(T)| \; |u|_{2,\infty,T} \,.
 \end{aligned}
\end{equation}
Here, $h(T)$ is the edge length of the square element $T$ and $\lesssim$ denotes the smaller or equal inequality up to a constant factor depending solely on the triangulation.
\end{corollary}

\textbf{Proof.} For the primal error weights we choose
$\tilde{u}_h = \Ih{2} u$, where $\Ih{2}$ denotes the piecewise bi-quadratic
Lagrangian interpolation.
Then standard interpolation estimates \cite[Theorem 3.1.5]{Ciarlet1978} give
\begin{align*}
 & \Vert u - \Ih{2} u \Vert_{0,2,T}          \lesssim h(T)           |u|_{1,2,T} \,, \\
 & \Vert u - \Ih{2} u \Vert_{0,2,\partial T} \lesssim h(T)^{\frac12} |u|_{1,2,\partial T} \,.
\end{align*}
Next, we estimate the error with respect to the control parameter $m$. The estimate for $\theta$ is completely analogous.
We choose $\tilde{m}_h = m[\efftensor \eps{\Ih{2} \strain}] =: \mathbf{m}(\lambda_1(\stress_h),\lambda_2(\stress_h))$,
i.e. $\tilde{m}_h$ is obtained by the pointwise evaluation of the interpolated solution $u$, the computation of the
corresponding stress $\stress_h=\efftensor \eps{\Ih{2} \strain}]$ and its eigenvalues and finally inserting them into formulae \ref{eq:explicitparameters}.
Then, $\sup_{x \in T} | \mathbf{m}(\lambda_1(\stress),\lambda_2(\stress))-$ $\mathbf{m}(\lambda_1(\stress_h),\lambda_2(\stress_h)) |$
is attained for some $\bar{x} \in T$ with stresses $\bar{\stress}=\stress(\bar x)$, $\bar{\stress}_h=\stress_h(\bar x)$.
From the Lipschitz continuity of the function $\mathbf{m}$ we deduce
\[
         \left|\mathbf{m}(\lambda_1(\bar{\stress}),\lambda_2(\bar{\stress}))-\mathbf{m}(\lambda_1(\bar{\stress}_h),\lambda_2(\bar{\stress}_h))\right|
\lesssim \left(|\lambda_1(\bar{\stress})-\lambda_1(\bar{\stress}_h)|^2+|\lambda_2(\bar{\stress}) - \lambda_2(\bar{\stress}_h)|^2 \right)^{\frac12} \!.
\]
Furthermore, the eigenvalues are Lipschitz continuous functions of the underlying matrices.
In particular the Wielandt-Hoffmann inequality \cite{HoWi53} for $p=2$ yields
\begin{align*}
&\left(|\lambda_1(\bar{\stress})-\lambda_1(\bar{\stress}_h)|^2+|\lambda_2( \bar{\stress}) - \lambda_2(\bar{\stress}_h)|^2 \right)^{\frac12}
\leq  \Vert \bar{\stress} - \bar{\stress}_h \Vert_F \,.
\end{align*}
The elasticity tensor $\efftensor$ represents a  bounded linear operator, thus
\[
 \left\Vert \stress - \stress_h \right\Vert_{0,\infty,T}
 = \left\Vert C (\eps{\strain} - \eps{\Ih{2} \strain} ) \right\Vert_{0,\infty,T}
 \leq ||| C ||| \;
 \Vert \eps{\strain} - \eps{\Ih{2} \strain} \Vert_{0,\infty,T}
\]
with $||| C |||$ representing the associated operator norm of the tensor $C$. Passing from the symmetrized strain tensor
to the full Jacobian yields
\[
 \Vert \eps{\strain} - \eps{\Ih{2} \strain} \Vert_{0,\infty,T} \lesssim |\strain - \Ih{2} \strain |_{1,\infty,T} \,.
\]
Finally, by standard $L^\infty$ estimates \cite[Theorem 3.3.7]{Ciarlet1978}, we achieve
\[
 |u - \Ih{2} u |_{1,\infty,T} \lesssim h(T) |\ln h(T)| \; |u|_{2,\infty,T} \,.
\]
\qed

\begin{remark}
 The presented argument would not work for a weighting term like $\left\Vert \alpha-\tilde{\alpha}_h \right\Vert_{0,\infty,T}$
 involving the rotation parameter as it depends on the eigenvectors of the local stress and the eigenvectors
 do not depend continuously on the matrix.
 In fact this motivates the elimination of the term involving the rotation parameter in Theorem~\ref{theorem:estimates}.
\end{remark}

\section{Numerical treatment} \label{sec:localapprox}
For a practical numerical scheme for the weighting terms $\omega_T$ that were a priori estimated in section~\ref{sec:apriori} we ask for an
approximation based on the numerical solution $u_h$.
One approach would be to exploit higher regularity of the solution by computing a second discrete solution using a higher order scheme.
However this is usually computationally demanding. To keep the cost for the evaluation of the error estimates low \cite{BeRa97} suggests to merely interpolate
the discrete solution to a higher order polynomial.
As we use bi-quadratic finite elements to overcome checkerboard instabilities, cf.~\cite{BoJo98}, we therefore construct polynomials of degree
four on patches of four neighboring elements from the local degrees of freedom of the discrete elastic displacement $u_h$.

In the discrete problem \eqref{eq:discrete} the coefficient functions $q$ and thus $\efftensor$ vary from point to point.
However it is known, see \cite[Theorem 4.1.6]{Ciarlet1978}, that the convergence order of the finite element scheme is unaltered
if the chosen quadrature order is high enough, which is the case in our computations.
To this end the parameters $\alpha$, $m$ can easily be evaluated at the quadrature points. The same holds true for the density parameter as well.
However it needs to be piecewise constant on each element to avoid
checkerboard instabilities, see \cite{BoJo98}. This is ensured by an element wise averaging of all
values obtained on quadrature points.
In detail we proceed as follows.
In case of the ratio parameter $m$ (and the rotation parameter $\alpha$ which however was eliminated from the error estimate)
the interpolated discrete deformation $u_h$ can be used again to locally obtain parameters $\alpha$ and $m$ for fixed $\theta$.
For the density parameter the piecewise constant approximation can be used to construct a bilinear profile on a patch of four neighboring
elements.  Let us summarize the approximation of the weights in the following remark:
\begin{remark}[Approximation of local weights]
To approximate the weights for the primal and the control error indicator on every element $T$ and to derive a suitable error indicator
we consider a patch $\mathcal{S}$ of the four neighboring elements including $T$ and proceed as follows:
\begin{itemize}
 \item[-] For the primal error indicator $5 \times 5$ degrees of freedom on each element of the patch are used to construct a bi-quartic
       polynomial on the patch.
       This interpolation $\Ih{4} \strain_h$ is then compared to the original approximation leading to
       \[
         \omega_T^{(u)}            \approx \left\Vert u_h-\Ih{4}u_h \right\Vert_{0,2,T}\,, \qquad
         \omega_{\partial T}^{(u)} \approx \left\Vert u_h-\Ih{4}u_h \right\Vert_{0,2,\partial T}\,.
       \]
 \item[-]
 For the control error indicator the ratio parameter $m$ is
 computed from the interpolated elastic solution above and compared to the original value;
 for the density parameter $\theta$ four piecewise constant values are used
 to construct a bilinear profile on the patch which is
 then compared to the original values. In explicit
       \[
         \omega_T^{(m)}      \approx \left\Vert m[u_h]   - m[\Ih{4}u_h]   \right\Vert_{0,\infty,T}\,, \qquad
         \omega_T^{(\theta)} \approx \left\Vert \theta_h - \Ih{1}\theta_h \right\Vert_{0,\infty,T}\,.
       \]
\end{itemize}
\end{remark}

\section{Implementation} \label{sec:impl}
For our computations we use an adaptive regular quadrilateral mesh provided by the \texttt{QuocMesh} library
\footnote{\url{http://numod.ins.uni-bonn.de/software/quocmesh}}.
Refinement is done via
uniform refinement of cells and handling of constrained hanging nodes. To obtain the optimal sequential lamination microstructure we
reimplemented the alternating algorithm suggested in \cite{AlBoFr97}.
The effective tensors are regularized by setting $C_{44} = 10^{-2}$. To ensure coercivity we restrict $m$ to the range
$[\varepsilon,1-\varepsilon]$ and $\theta$ to $[\varepsilon,1]$ with $\varepsilon=10^{-3}$.
The algorithm is terminated once the change in the compliance cost value is below $10^{-7}$.
As already mentioned we use bi-quadratic elements to overcome checkerboard instabilities.
For the numerical integration it turned out to be sufficient to use a Gauss quadrature rule of order 5.
The lamination parameters for the interpolated solution cannot be directly computed from the pointwise stress
$\efftensor \, \eps{\Ih{4} \strain_h}$ as it in turn depends on the unknown parameters.
Here we use a Newton method to solve the equation
\[
 \bar{\etensor}[m(\lambda_i),\bar{\theta}] \, \eps{\Ih{4} \strain_h} - Q(\alpha) \begin{pmatrix} \lambda_1 & \\ & \lambda_2 \end{pmatrix} Q(\alpha)^\top = 0
\]
for the rotation parameter $\alpha$ and the eigenvalues $\lambda_i$ of the stress (from which $m$ can be computed).
The density $\bar{\theta}$ is fixed to the original value found at the current point.

For the adaptive scheme we use a D\"orfler marking strategy \cite{Do96} with a fraction of $40\%$ of the total error to be associated with
the set of elements to be refined.

\section{Numerical results} \label{sec:results}
\begin{figure}[!ht]
\includegraphics[width=.18\linewidth]{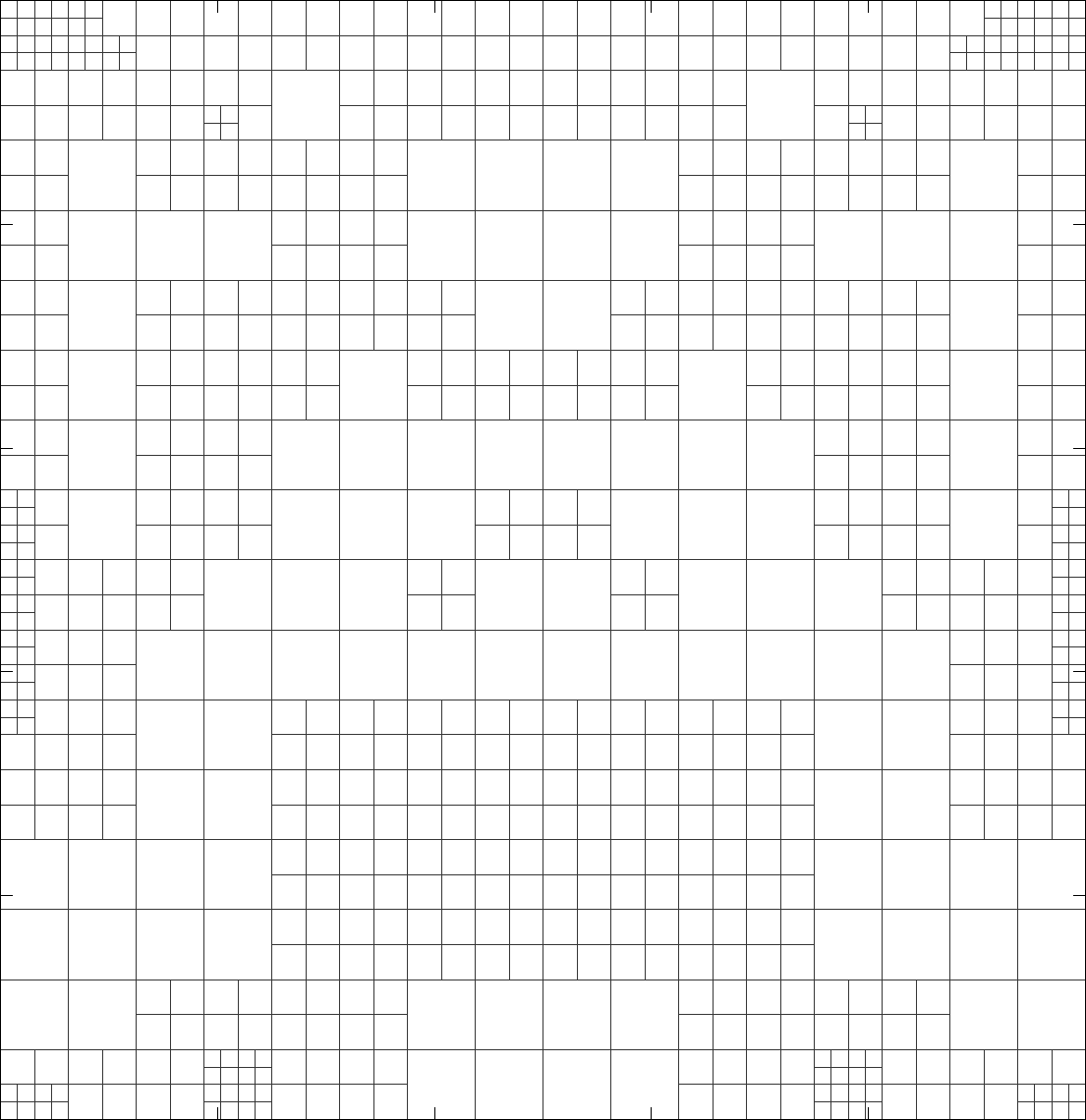}
\hfill
\includegraphics[width=.18\linewidth]{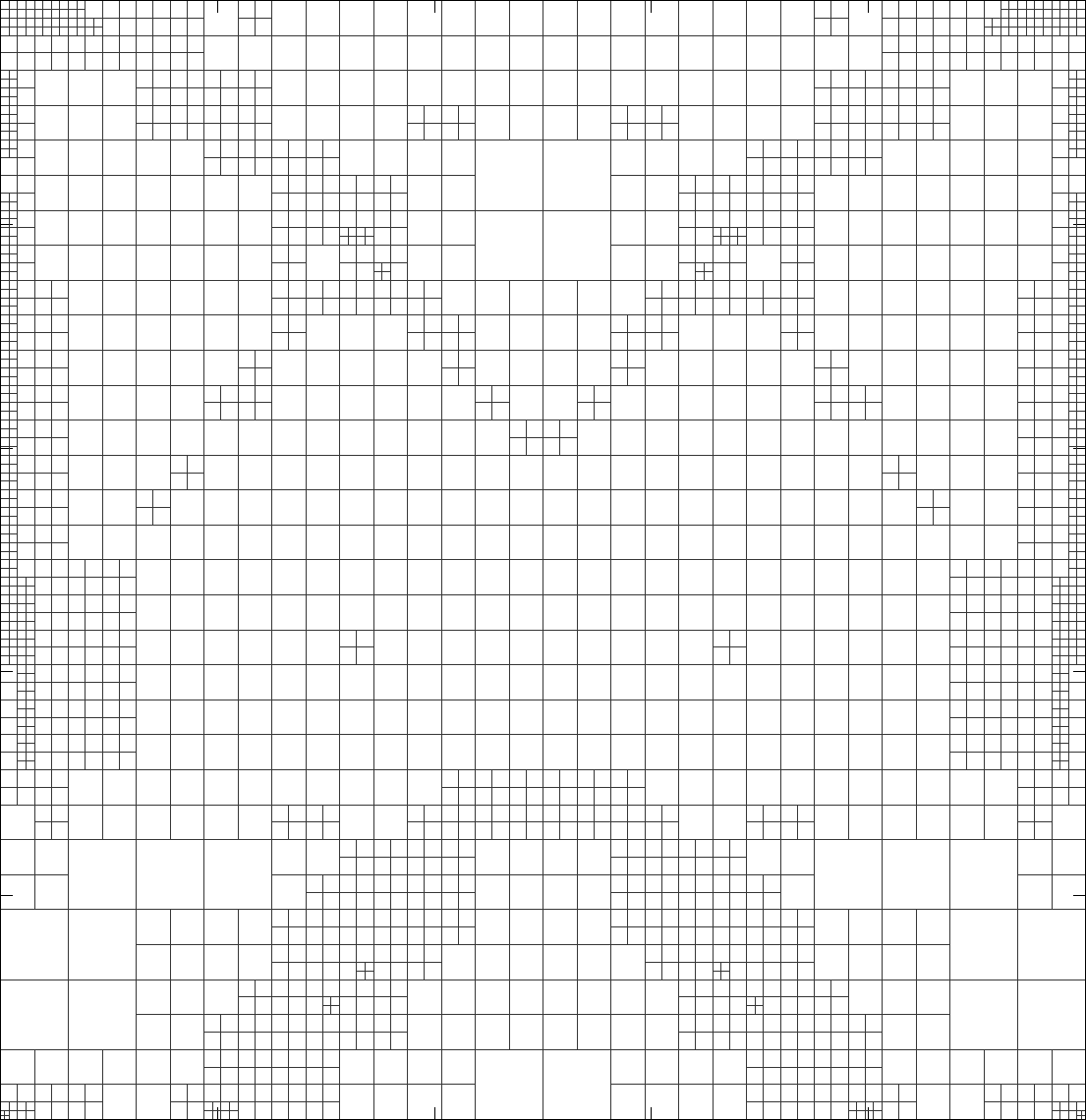}
\hfill
\includegraphics[width=.18\linewidth]{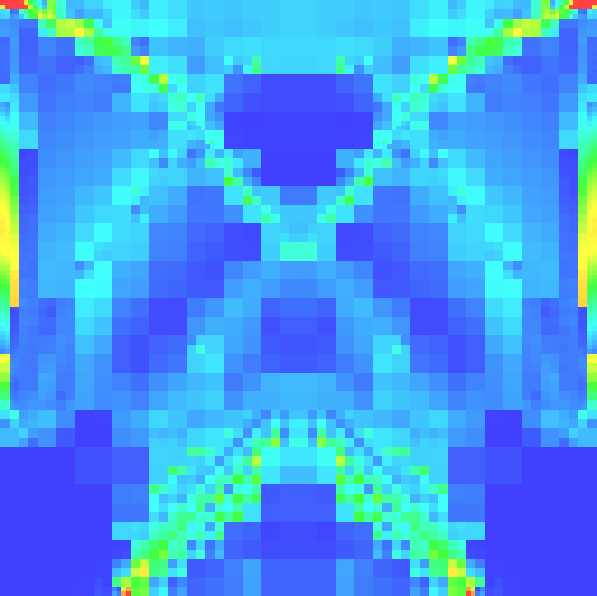}
\hfill
\includegraphics[width=.18\linewidth]{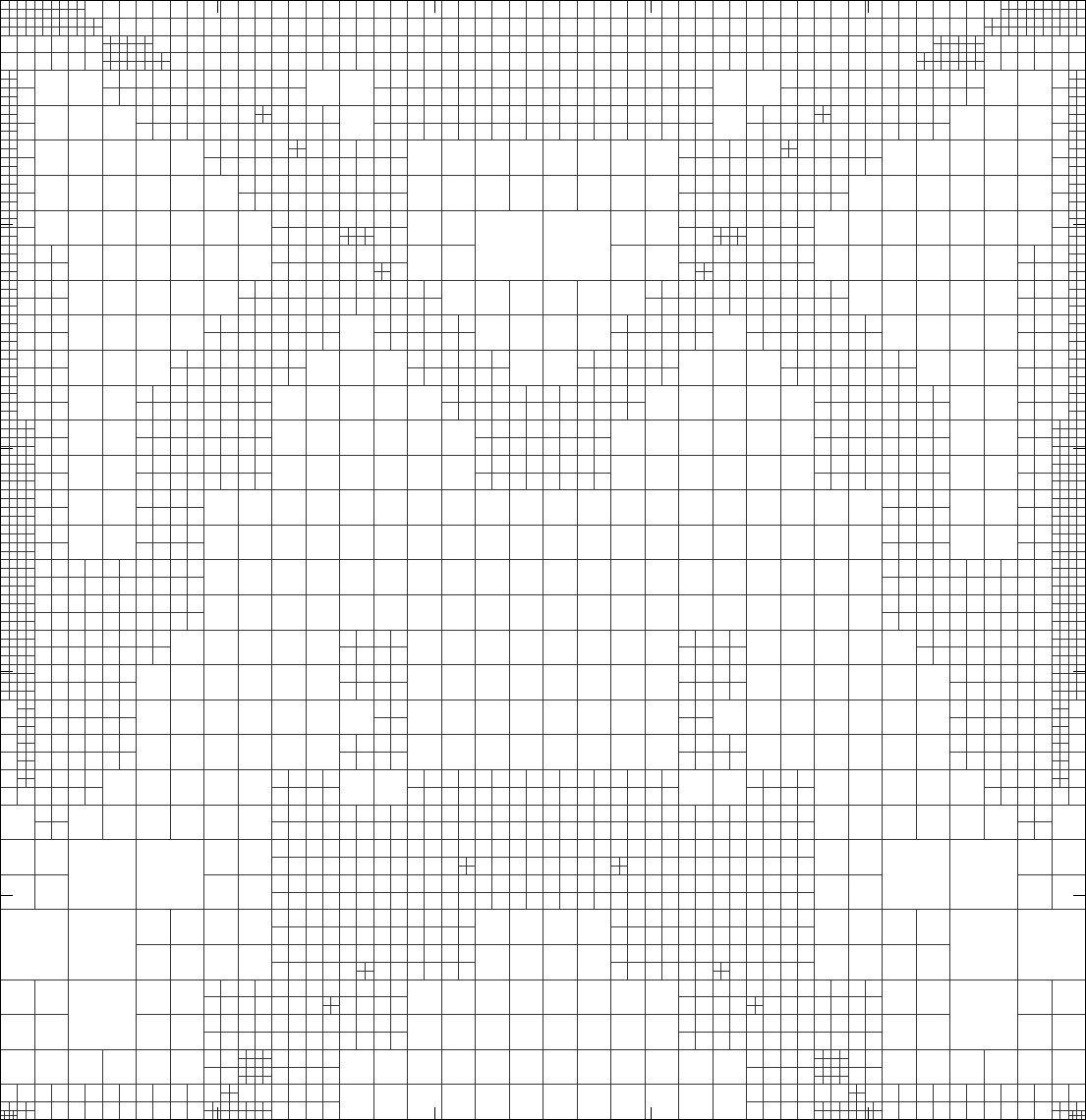}
\hfill
\includegraphics[width=.18\linewidth]{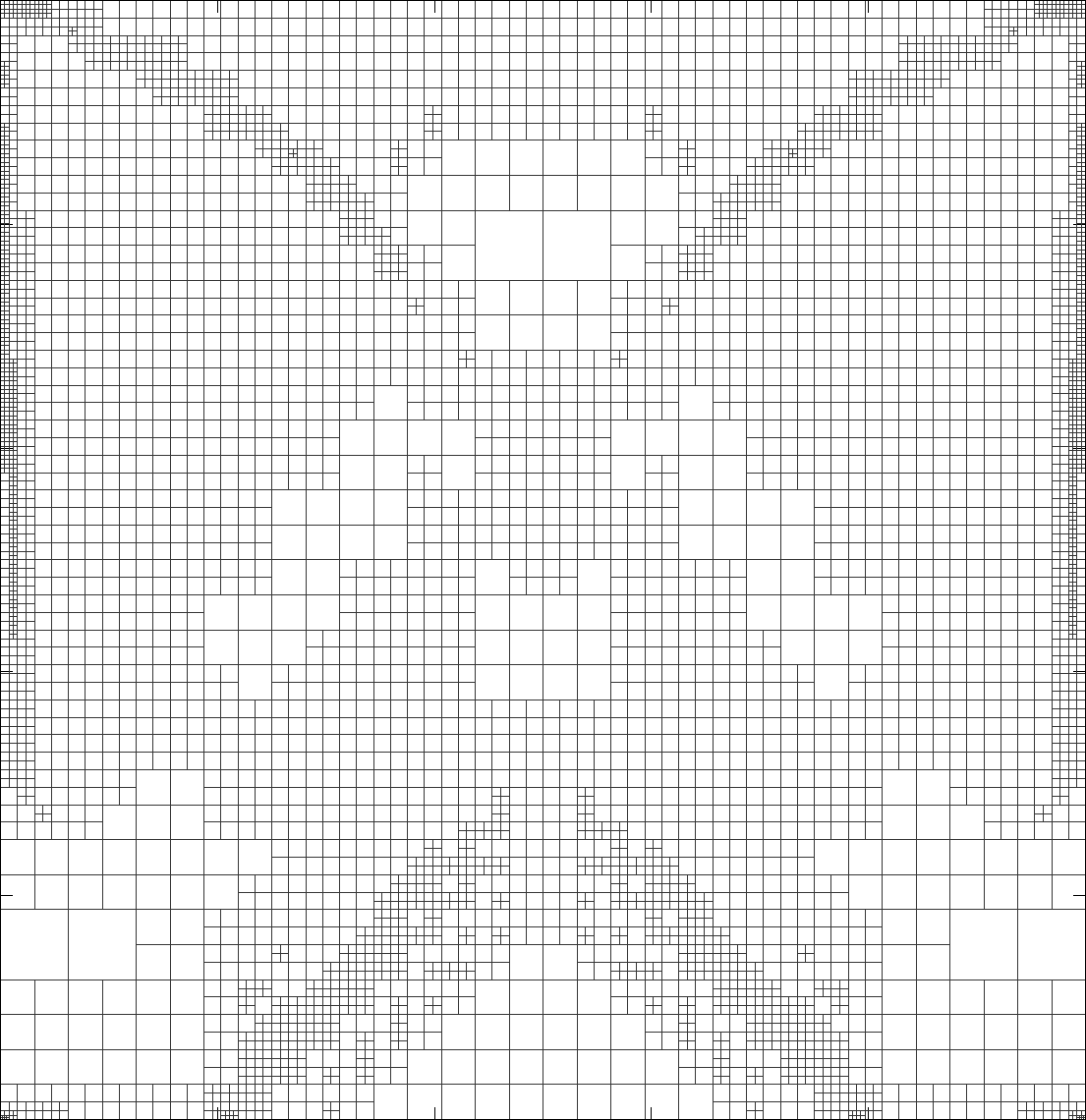}\\[3pt]
\includegraphics[width=.18\linewidth]{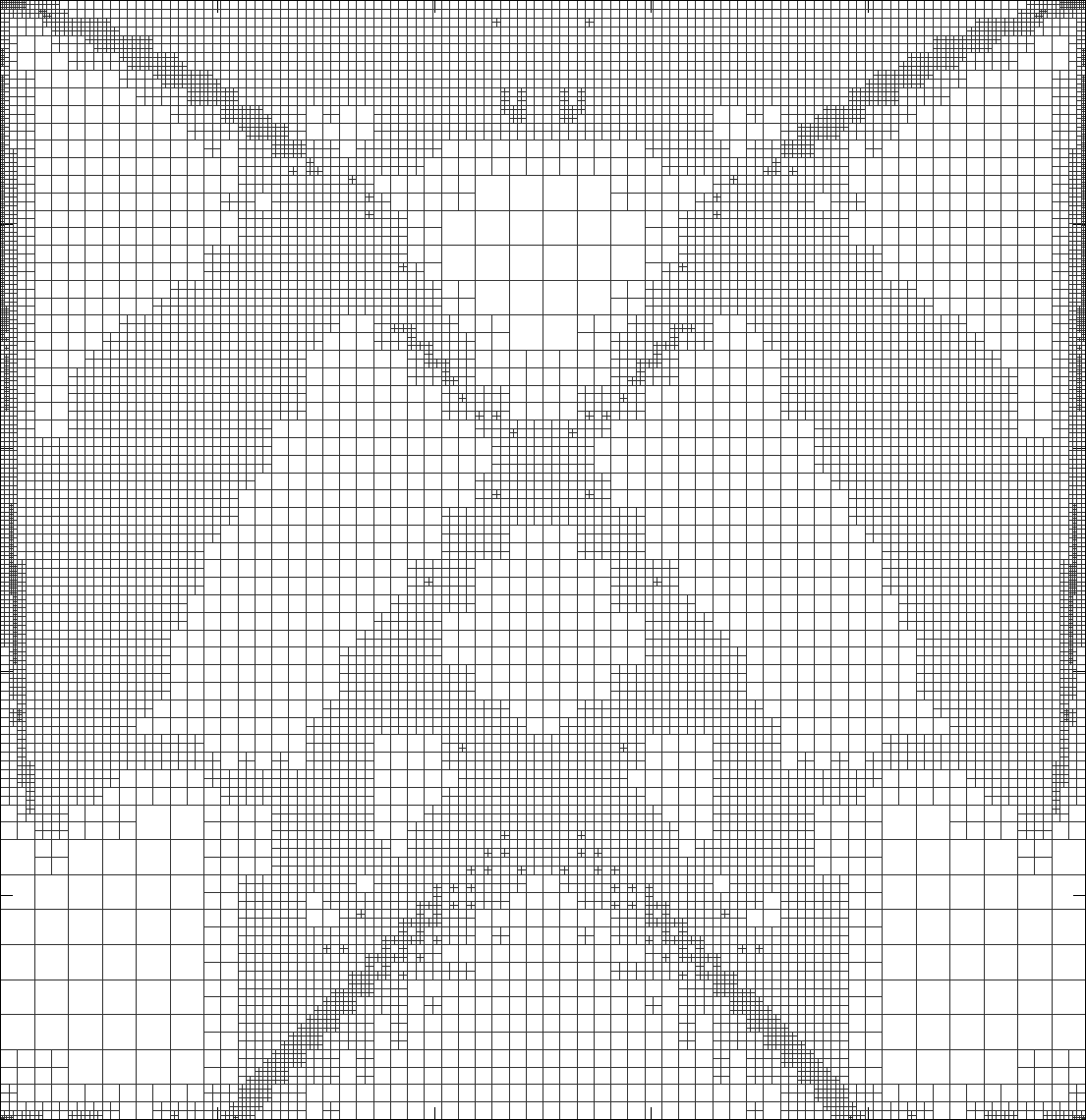}
\hfill
\includegraphics[width=.18\linewidth]{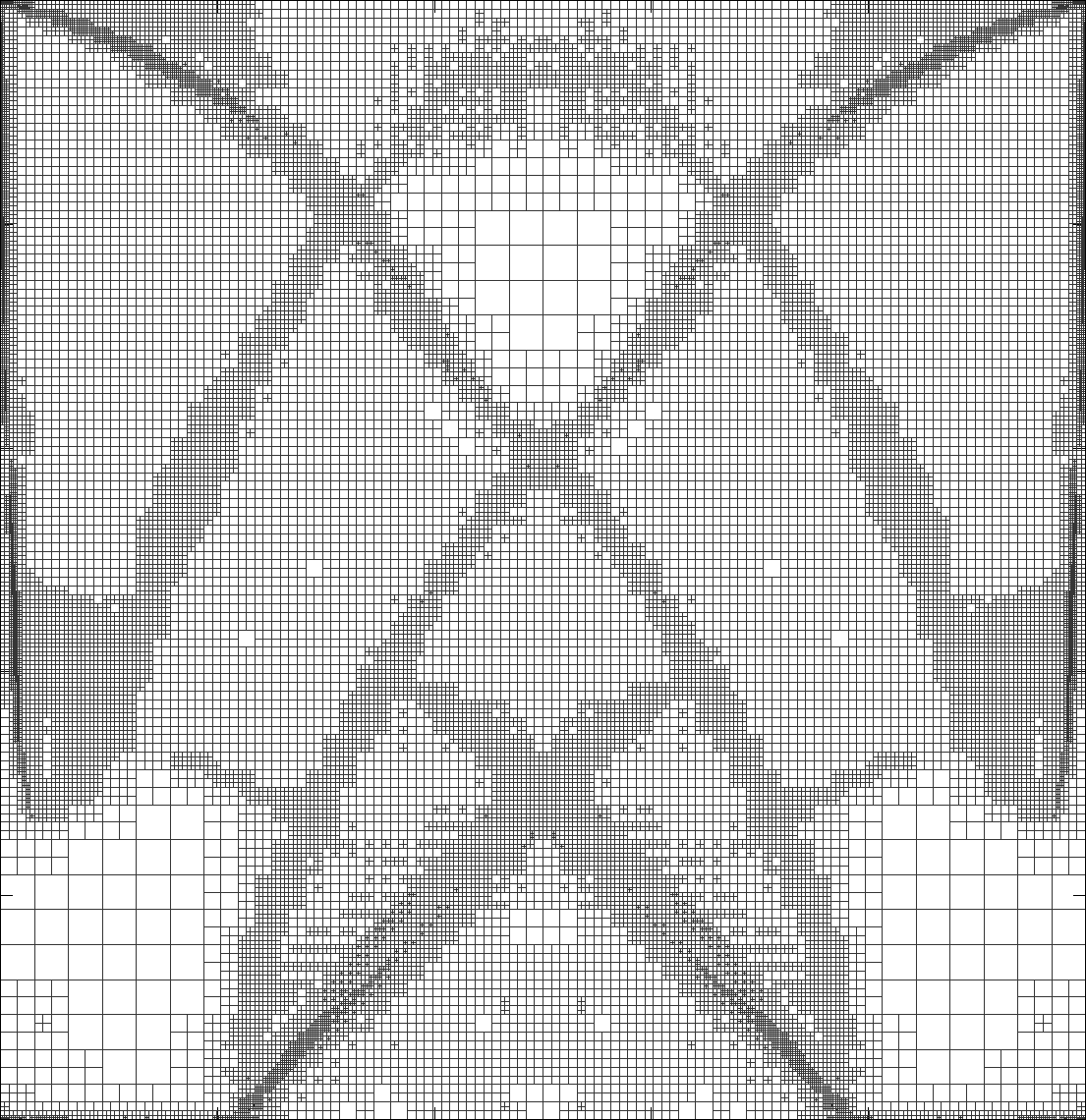}
\hfill
\includegraphics[width=.18\linewidth]{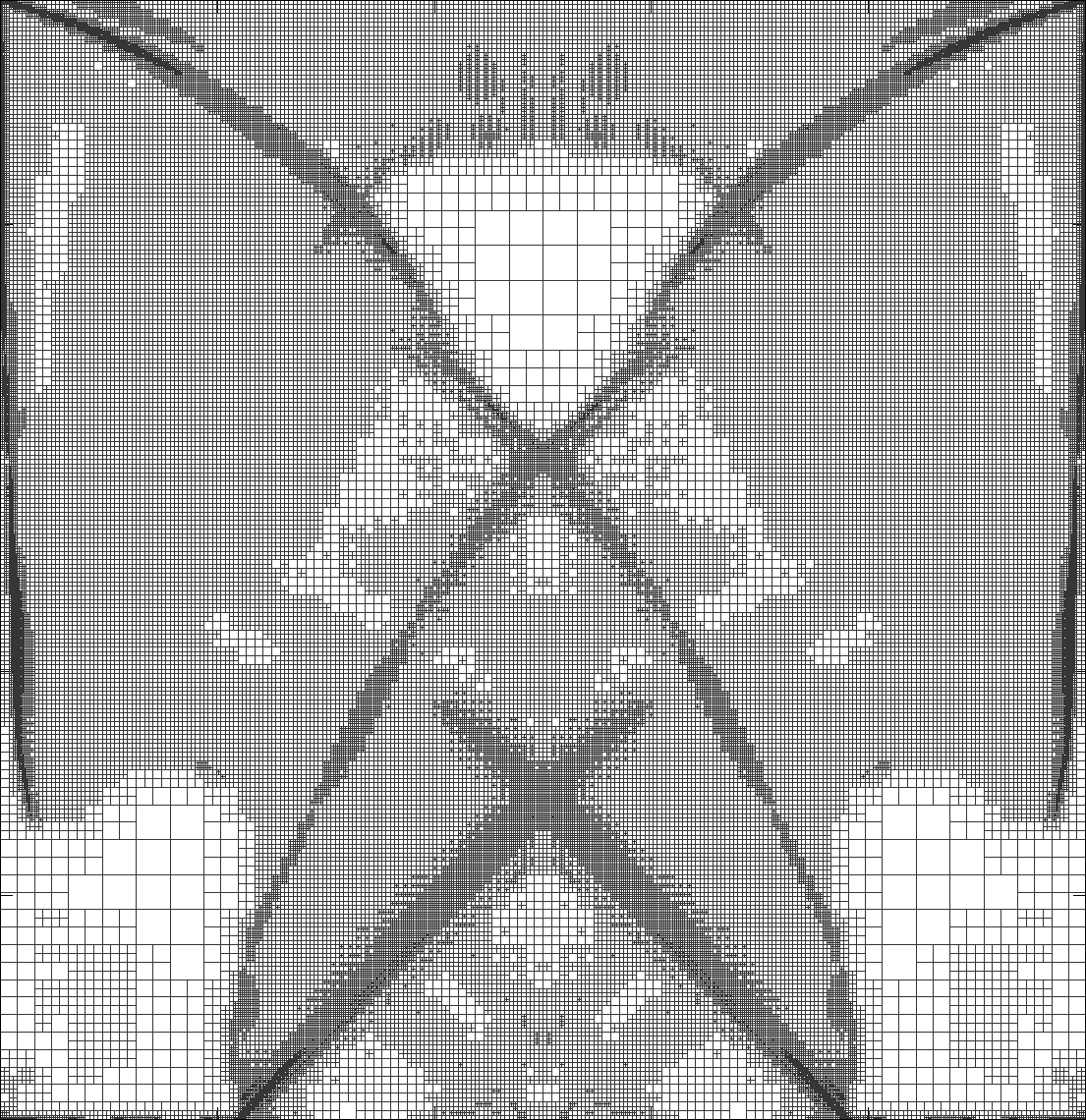}
\hfill
\includegraphics[width=.18\linewidth]{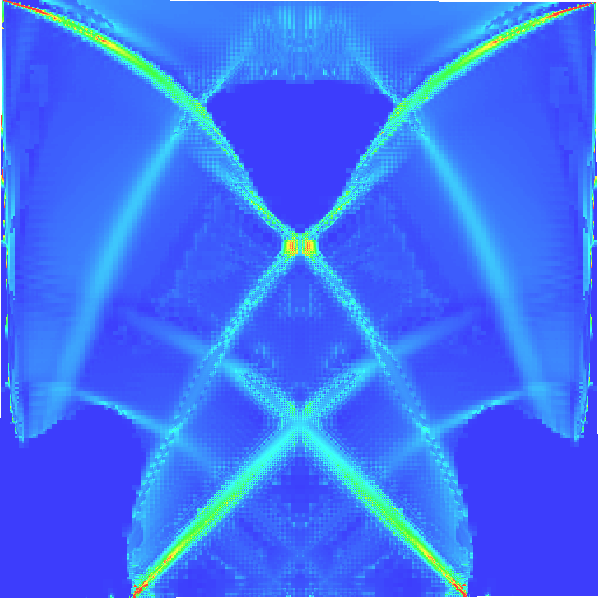}
\hfill
\includegraphics[width=.18\linewidth]{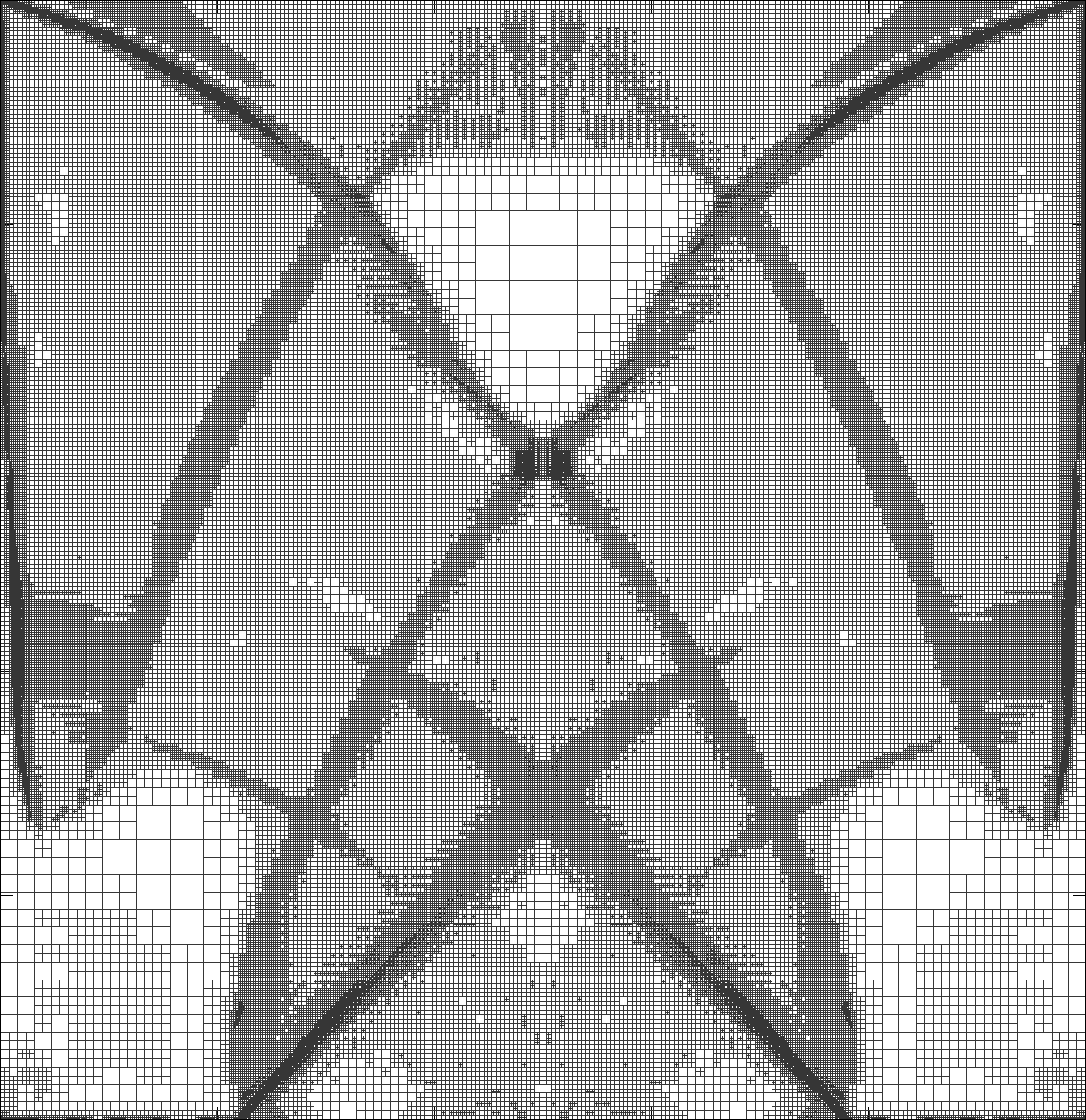}
\caption{Adaptive meshes after 4, 7, 8, 10, 13, 16, 19, and 20 refinement steps obtained for the carrier plate scenario.
Color plots show the error indicator after 7 and 19 steps, leading to the subsequent grid.}
\label{fig:carrier}
\end{figure}

As the first example we consider a {\it carrier plate under shearing} computed on the unit square., cf.~Fig.~\ref{fig:carrier}.
The volume constraint was set to $33\%$.
First we ran the alternating algorithm until convergence on uniform meshes of level $l \in \{2,\ldots,10\}$.
The obtained values can be used to extrapolate the asymptotic value of the compliance cost by assuming the following
dependence on the grid width:
\[
 J[u_h] = J^\ast + c \, h^p \,.
\]
Using a least squares fit we found the parameters $J^\ast = 1{.}8399$, $c = 1{.}7645$, and $p = 1{.}0484$. 
$J^\ast$ is used to replace the exact cost value when assessing error reduction in the following.

\begin{figure}[!ht]
\hfill
\includegraphics[width=.4\linewidth]{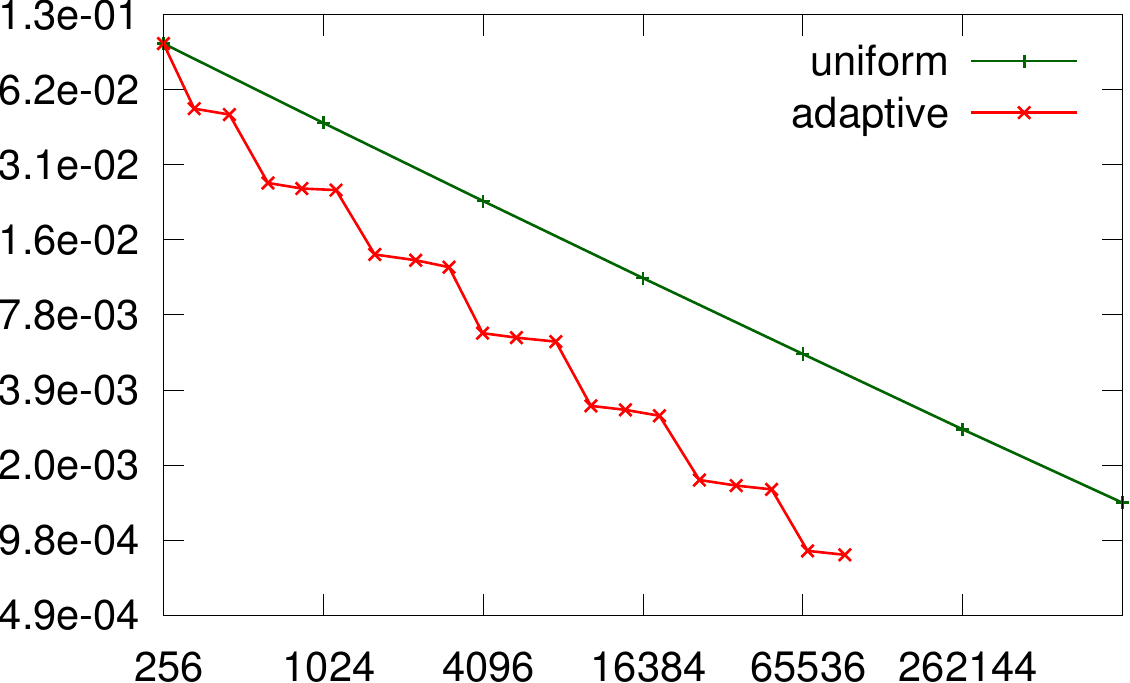}
\hfill\mbox{}
\caption{
Difference to extrapolated value $J^\ast$ plotted
over the number of elements, once for a uniform refinement (green) and the adaptive refinement strategy (red).
        }
\label{fig:carrier2}
\end{figure}

We ran the discussed adaptive scheme for the carrier plate scenario, see Figures~\ref{fig:carrier} and \ref{fig:carrier2}.


\begin{figure}[!ht]
\includegraphics[width=.3\linewidth]{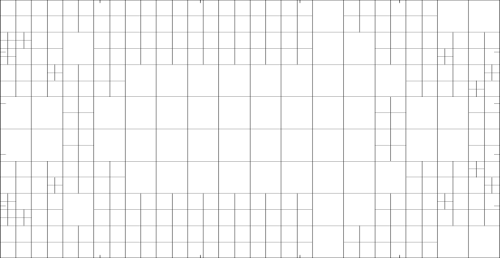}
\hfill
\includegraphics[width=.3\linewidth]{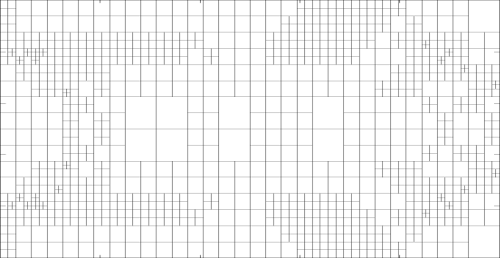}
\hfill
\includegraphics[width=.3\linewidth]{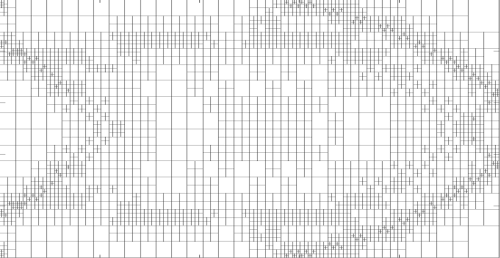}\\[3pt]

\includegraphics[width=.3\linewidth]{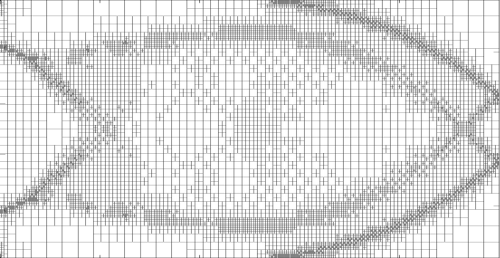}
\hfill
\includegraphics[width=.3\linewidth]{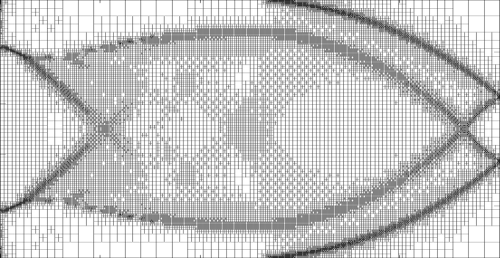}
\hfill
\includegraphics[width=.3\linewidth]{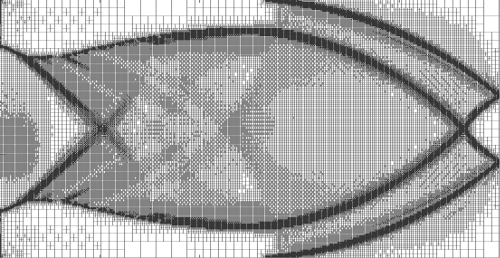}\\
\mbox{\rule{\linewidth}{.4pt}}\\[12pt]

\includegraphics[width=.3\linewidth]{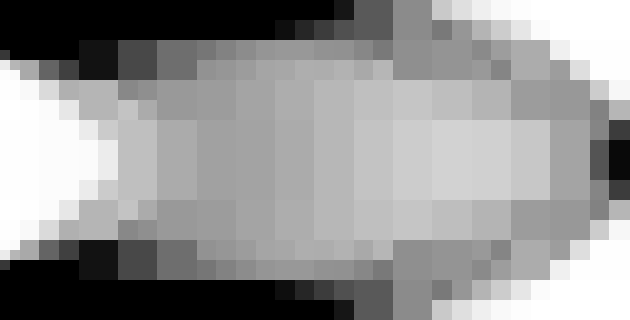}
\hfill
\includegraphics[width=.3\linewidth]{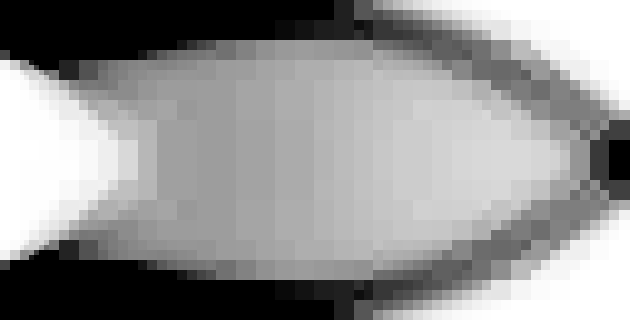}
\hfill
\includegraphics[width=.3\linewidth]{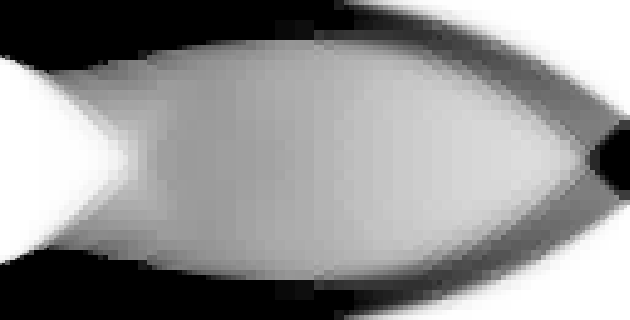}\\[3pt]

\includegraphics[width=.3\linewidth]{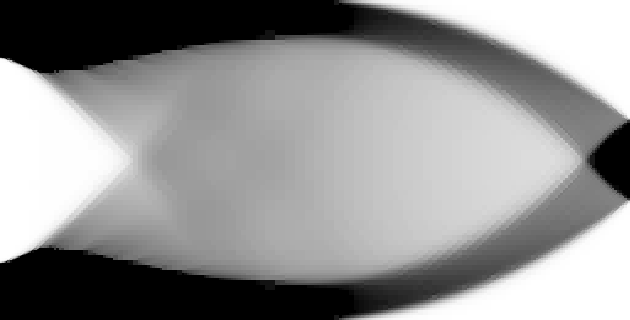}
\hfill
\includegraphics[width=.3\linewidth]{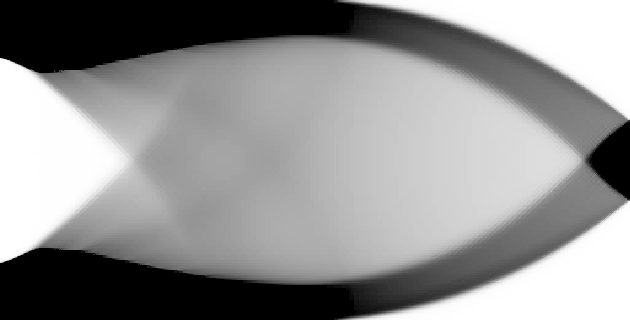}
\hfill
\includegraphics[width=.3\linewidth]{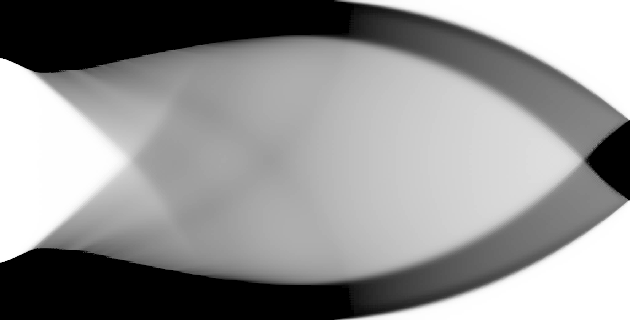}\\
\caption{At different refinement steps of the adaptive method applied to the cantilever example the current mesh (top) and the corresponding density $\theta$ (bottom) are displayed (from left to right: step 4, 8, 12, 16, 20, and 24).}
\label{fig:cantilever}
\end{figure}

Next, we take into account a {\it cantilever scenario}. Here we consider the domain $D = [0,2] \times [0,1]$
with  Dirichlet boundary condition on the left hand side and a downward pointing force in the middle on the right hand side. The
volume constraint is set to $50\%$, cf. Figure~\ref{fig:cantilever}.

\begin{figure}[!hb]
\hfill
\includegraphics[width=.4\linewidth]{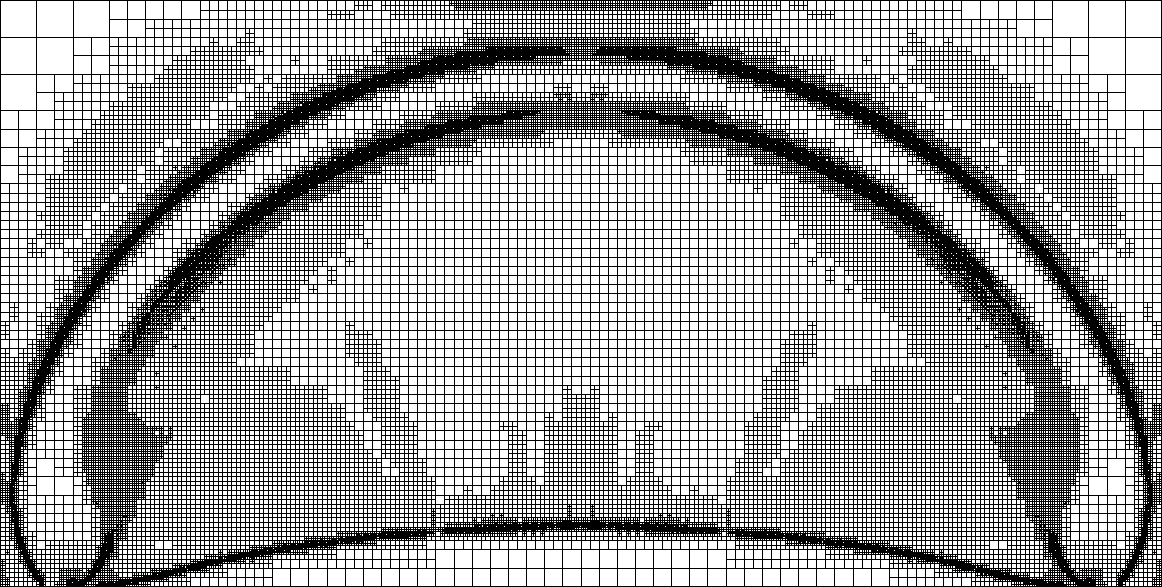}
\hfill
\includegraphics[width=.4\linewidth]{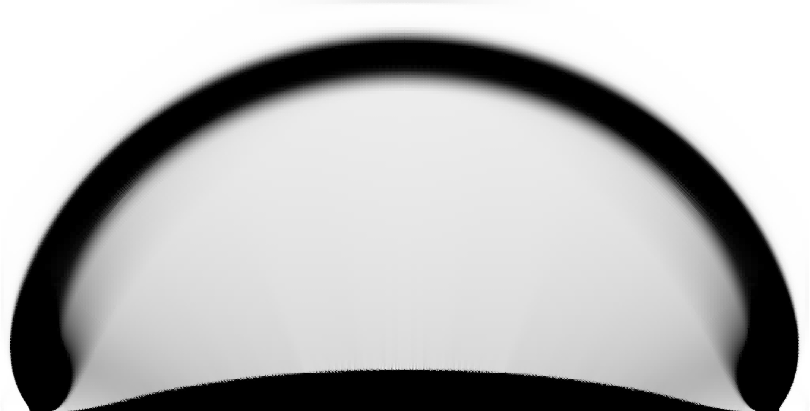}
\hfill\mbox{}
\caption{Adaptively refined grid and corresponding density $\theta$ for the bridge scenario.}
\label{fig:bridge}
\end{figure}

As a third scenario we investigate a {\it bridge} within the domain $D = [0,2] \times [0,1]$ in Figure~\ref{fig:bridge}.
We prescribe roller boundary conditions on short strips at the lower boundary in the left and right corners, \ie\
the displacement in $y$-direction is enforced to be $0$, compare for instance \cite{Al02}. 
Additionally the node in the lower left corner is subject
to full homogeneous Dirichlet boundary conditions to ensure a unique solution. The lower boundary in between is loaded
in vertical direction. The volume constraint is set to $33\%$.

\section{Conclusion} \label{sec:sum}
In this paper we have applied the dual weighted residual approach to an elastic multi scale shape optimization
problem using sequential laminates as an underlying optimal microstructure. The microstructure
enters the PDE via effective material properties whose explicit formulae allow to compute
sensitivities of the cost functional with respect to the describing parameters of the microstructure.
The derived goal-oriented error estimate enables to appropriately steer the adaptive refinement
leading to a substantial reduction of the number of degrees of freedom  compared to uniform meshes and the same accuracy.
In particular the goal oriented error estimation outperforms a residual type estimator solely for the elastic displacement problem.
The main ingredient  is the error term incorporating the sensitivity of the elastic energy density \wrt\ to the local material density.
This is in accordance with the observation that the rigidity of a shape is mostly improved by an optimal distribution of the available material.
The adaptive scheme leads to sharply resolved interfaces, separating regions with
almost no material from regions with bulk material or an intermediate material density.\\

\section*{Acknowledgments}
This work has been supported by the \emph{Deutsche Forschungsgemeinschaft} through
\emph{Collaborative Research Centre 1060 -- The Mathematics of Emergent Effects}.

\bibliographystyle{alphadin}
\bibliography{bibtex/all,bibtex/own,bibtex/library,additional}

\end{document}